\documentclass[12pt]{article}
\usepackage[english]{babel}
\usepackage{amsmath,amsthm}
\usepackage{amsfonts}
\usepackage{authblk}
\usepackage{hyperref}
\usepackage[lined,boxed,linesnumbered]{algorithm2e}
\usepackage{tikz}
\usepackage{graphicx}
\usepackage{subcaption}
\usepackage{multirow}
\usepackage{placeins}
\usepackage{pdflscape}
\usepackage{longtable}
\usepackage[utf8]{inputenc}
\usepackage{comment}
\usepackage{geometry}
\newtheorem{thm}{Theorem}[section]

\theoremstyle{definition}

\theoremstyle{remark}
\newtheorem{rem}[thm]{Remark}
\numberwithin{equation}{section}

\usepackage{color}

\newcommand{\BI}{\begin{itemize}}
\newcommand{\EI}{\end{itemize}}

\newcommand{\st}{\mathop{\rm s.t.}}

\newcommand{\crosse}{\gamma}				

\newcommand{\Real}{\mathbb{R}}

\newcommand{\hc}{\hat{\beta}}

\newcommand{\uset}{\mathcal{U}}
\newcommand{\ruset}{\mathcal{RU}}

\newcommand{\exclude}[1]{}

\title{Robust Price Optimization of Multiple Products under Interval Uncertainties}
\author{Mahdi Hamzeei\footnote{hamzeei.m@gmail.com} \hspace{0.25in}
Alvin Lim\footnote{alvin.lim@nielseniq.com, corresponding author} \hspace{0.25in}
Jiefeng Xu\footnote{jiefeng.xu@nielseniq.com, principal corresponding author} }

\affil{Precima, a NielsenIQ Company\\ 200 West Jackson, Blvd., Chicago, IL 60606}

\date{}

\begin{document}

\maketitle


\begin{abstract}
In this paper, we solve the multiple product price optimization problem under interval uncertainties of the price sensitivity parameters in the demand function.
The objective of the price optimization problem is to maximize the overall revenue of the firm where
the decision variables are the prices of the products supplied by the firm.
We propose an approach that yields optimal solutions under different variations
of the estimated price sensitivity parameters. We adopt a robust optimization approach by building
a data-driven uncertainty set for the parameters, and then construct a deterministic
counterpart for the robust optimization model.
The numerical results show that two objectives are fulfilled: the method reflects the uncertainty
embedded in parameter estimations, and also an interval is obtained for optimal prices.
We also conducted a simulation study to which we compared the results of our approach.
The comparisons show that although robust optimization is deemed to be conservative,
the results of the proposed approach show little loss compared to those from the simulation.

\end{abstract}


\section{Introduction}
\label{sec:intro}

In today's marketplace, price is considered as a key driver for optimizing revenue as it significantly affects the demand on products. Thus, the price optimization problem has long been studied as a significant decision-making problem for companies. There are numerous books, e.g., Phillips \cite{PRO:RoberPhillips:2005}, Talluri and van Ryzin \cite{Springer:talluri:05} and {\" O}zer and
Phillips \cite{Ox:ozer:12}, that provide more comprehensive reviews on various price optimization models.

In this paper, we present a data-driven {\it robust optimization} (RO) based approach for the multi-product price optimization problem (MPPO) for a food supplier. The objective of the problem is to maximize the revenue of selling the products supplied by the firm while
the decision variables are the prices of these products. We assume that for each product, a response function describes the demand for the product that depends on the price of the product itself as well as on the corresponding prices of its complementary and substitute products.
The parameters of the demand function, which reflect the consumers' demand sensitivity to prices are typically estimated through various statistical approaches
with associated confidence intervals. In practice, the price optimization problem is solved by assuming that the mid-points of the confidence intervals as the point estimators of the parameters to compute an approximation of the demands and the associated revenues at various prices. However, the solution obtained by this approximate
problem may not reflect well the uncertainty involved in the parameter estimates.

In most of the existing studies, the demand functions with price variables are built with parameters whose values are obtained either from a priori knowledge or from exact estimates using existing data. Nevertheless, in practice the estimates of these parameters may take significantly different values in different circumstances, i.e., different customers may have different sensitivities to price, or customers may have different price reactions at different periods of the year. By neglecting such uncertainties, the demand functions may result in potentially inaccurate or even incorrect demand estimates at various prices. As suggested by a few researchers, e.g., Nahmias~\cite{MG:Nahmias:05}, Simchi-Levi et.~al.~\cite{MG:Simchi:04} and Sheffi~\cite{MIT:Sheffi:05}, many point estimates used in modern operations management decision models are wrong or meaningless. Instead, they recommend that the point estimates be replaced by their corresponding ranges to handle uncertainties.

Although in the field of revenue management, most conventional approaches aim to maximize the expected value of the revenue (e.g. Agrawal and Seshadri~\cite{Agrawal:MSOM:00}), the risk-neutrality of the expected value may not well capture risk preferences of decision makers in many applications. Therefore, Levin et.~al.~\cite{Levin:OR:08}, among the others, proposed the use of value-at-risk constraints to handle uncertainty in a more risk-averse fashion. However, one needs to estimate a probability distribution for the uncertain value and this often becomes a challenging task in practice. Furthermore, these parameters with probability information impose additional complexities to the underlying optimization models.

Our goal is to propose an approach to deal with the uncertainty involved in the estimated parameters of the demand functions. Our RO approach replaces point estimates by interval forecasts and requires no specific probability distribution. Instead, we choose a hyperplane uncertainty set based on constraining the deviation of uncertainty parameters from the nominal values.  As an additional benefit, our RO approach addresses the risk-averse tendencies of decision makers.

To the best of our knowledge, the idea of RO was first introduced by Soyster \cite{Soyster:OR:73}. Since then, especially with the advent of efficient algorithms to solve optimization problems, the field has attracted much more attention, and expanded significantly over the last few decades. RO addresses optimization problems with uncertainty in which the uncertainty model is not stochastic, but rather deterministic and set-based. For example, in 1995, Mulvey et.~al.~\cite{Mulvey:OR:1995} developed a robust optimization approach to solve large-scale systems for a set of decision-making problems. The paper introduces a general framework for achieving a solution that is robust in terms of both
feasibility and optimality with respect to {\it all} realizations of uncertain data. It also presents several applications solved using this approach.

The literature on computationally tractable RO models is rich. For instance, linear programming with ellipsoidal uncertainty sets has been addressed in \cite{Ben-Tal:1998:MOR}, \cite{Ben-tal:99:ORL}, \cite{Ben-Tal:00:MP}, \cite{ElGhaoui:siam:97} and \cite{Ghaoui:SIAM:98}. In such a case, the optimization problem will correspond to a set of conic quadratic programming optimization problems. Alternatively, RO models with polyhedral uncertainty sets which can be formulated with linear/integer variables are also well-studied (see  \cite{Ben-tal:99:ORL,Bertsimas:mp:03,Bertsimas:OR:04} for examples). For a thorough review of RO, we refer the interested readers to \cite{Bertsimas:SiamRev:11}.

One of the most crucial steps in applying RO approaches is to construct an uncertainty set such that it contains all, or almost surely all, possible variations of uncertain data. Bound constraints and ellipsoidal uncertainty sets have been commonly used in many treatments \cite{Bertsimas:ORL:04}. However, these approaches are often ad-hoc, with emphases on sets that preserve computational tractability. Bertsimas and Brown \cite{Bertsimas:OR:2009} and Bertsimas et.~al. \cite{Bertsimas:MP:2018} have proposed methods to construct polyhedral uncertainty sets using statistical hypothesis tests.

Previously, RO has been applied to several variations of the price optimization problem. Thiele~\cite{thiele:OO:06} applied an RO approach to the pricing problem of a single product over a finite time horizon with a capacity limit. The uncertainty in the paper is the demand of
the product which is assumed to be a function of price. The problem considered is different from the multiple products price optimization problem in our paper. In another paper, Thiele~\cite{thiele:jrpm:09} proposed an RO approach for a problem with multiple products where a capacitated resource constraint is enforced, and the uncertainty is on the demand of the products. In contrast to the uncertainty set in both papers, our paper directly addresses the uncertainty on the sensitivity to the underlying prices in the demand functions which auguments the uncertainty in the demand of the products.

In addition, Tien and Jaillet~\cite{Mai:arxive:19} studied an interesting pricing problem
with general extreme value (GEV) choice models for customers. They assumed that the
parameters of the choice model lie in an uncertainty set. They first consider
an unconstrained pricing problem, then present an alternative version, and analyze the pricing problem with over-expected-revenue-penalties. For the Multinomial Logit (MNL) demand function and a rectangular uncertainty set, the RO problem can be converted to a deterministic one that can be solved efficiently.

The rest of our paper is organized as follows. In Section~\ref{sec:method}, we explain the multi-product price optimization model. The underlying uncertainty as well as the approach we use to handle the uncertainty are explained
in Section~\ref{sec:uncertainty}. We then present our numerical experiments in Section~\ref{sec:NumExam}.
Finally, we conclude the paper with several findings  and point out some potential future research directions in
Section~\ref{sec:conclusion}.

\section{The Price Optimization Model}
\label{sec:method}

In this section, we first formulate the deterministic case of the multi-product price optimization problem which is modeled as a nonconvex quadratic problem with bilinear terms. We then present a relaxation of the problem into a linear program.

\subsection{The Deterministic MPPO Model}
For our MPPO, we assume a {\it linear price-response function}
as suggested by Phillips \cite{PRO:RoberPhillips:2005}:

\begin{equation}
\label{eq:baselinear}
d(p) = \alpha + \beta p
\end{equation}
where $\alpha>0$ and $\beta <0$ are the respective intercept and slope of the linear price-response function $d(p)$ which represents the sales volume as a function of price. The {\it satiating price} $P$, defined as the maximum allowed price at which demand
drops to zero, is given by $P = -\alpha / \beta $.
Phillips \cite{PRO:RoberPhillips:2005} points out that the {\it elasticity} of
 the linear price-response function is -$\beta p / (\alpha + \beta p)$ which ranges from
 0 at $p=0$ and approaches infinity as $p$ approaches $P$.

In reality, it is well known that the demand of a product is often influenced by not only the price of itself, but also the prices of its complementary and substitute products. Therefore we extend the linear function \eqref{eq:baselinear} to incorporate
the prices of complementary and substitute products. Let $I$ represent the set of product indexes.
The extended demand function for product $i$ can be expressed by introducing additional price variables along with complementary/substitute products:

\begin{equation}
\label{eq:genlinear}
d_i(p) = \alpha_{i} + \beta_{i}p_i + \sum_{j\in C_i}{\crosse_{ij} p_j}
\end{equation}
where $p_j$ is the price of the complementary/substitute product $j$, $\crosse_{ij}$ represents the corresponding cross-effect
of product $j$ on own product $i$, and $C_i$ is the set of associated complementary/substitute products for product $i$.

The MPPO is to maximize the total revenue subject to a set of constraints:

\begin{equation}
\label{prob:orig}
	\begin{array}{rl}
		\max\limits_p & \sum\limits_{i\in I}{p_i \cdot  (\alpha_{i} + \beta_{i}p_i + \sum\limits_{j\in C_i}{\crosse_{ij} p_j} ) }    \vspace{1mm}\\
		\st & Ap \leq b, \vspace{1mm} \\
		& l \leq p \leq u,
	\end{array}
\end{equation}
where $p = (p_1, p_2, \ldots, p_n) \in \Real^n$ (hence $n=|I|$) , $A\in\Real^{m\times n}$, $b\in \Real^m$,
and $l,u\in\Real_+^n$ are finite.  The constraints in \eqref{prob:orig} defined by matrix $A$ and vector 
$b$ represent business rules in practice.

We assume that the above problem is well-defined and {\it feasible}.  We also assume that for $i \in I$, $\beta_i < 0$, so the demand of a product decreases when its price increases, and $\gamma_{ij} < 0$ for any complementary product $j \in C_i$ and $\gamma_{ij} > 0$ for any substitute product $j \in C_i$.

Thus, the objective function
\[
f(p; \beta) = \sum_{i\in I}{\left( \alpha_{i} p_i  +  \beta_{i}p_i^2 +
	\sum_{j\in C_i}{\crosse_{ij} p_i p_j} \right) },
\]
is, in general, a nonconvex quadratic function with bilinear terms.

\subsection{A Relaxation of the Deterministic MPPO Model}
For the quadratic and bilinear terms, $p_i^2$ and $p_i p_j$, in \eqref{prob:orig}, we introduce two new sets of
variables $x_i=-p_i^2$ and $y_{ij}=p_i\cdot p_j$, respectively. With this substitution, (\ref{prob:orig}) is transformed into:

\begin{subequations}
\label{prob:reform1}
\begin{align}
		\label{prob:reform1:obj} \max\limits_{p,x,y} \quad & \sum\limits_{i\in I}{ \left( \alpha_{i} p_i - \beta_{i} x_i +
	\sum_{j\in C_i}{\crosse_{ij} y_{ij}} \right) } \vspace{1mm} & \\
		\label{prob:reform1:const1} \st \quad & Ap \leq b, & \vspace{1mm}  \\
		\label{prob:reform1:const2} & x_i \leq -p_i^2, & i\in I,  \\
		\label{prob:reform1:const3} & y_{ij} = p_ip_j, & i\in I, j \in C_i,  \\
		& l \leq p \leq u.
\end{align}
\end{subequations}
where the objective function \eqref{prob:reform1:obj} is now linear, but the constraints \eqref{prob:reform1:const2} and \eqref{prob:reform1:const3} contain quadratic and bilinear terms, respectively.

Notice that for a given $i \in I$, since $\beta_i < 0$ and $p_i > 0$, $x_i < 0$ and $-\beta_i x_i < 0$.  Because we are dealing with a maximization problem, it follows that $x_i$ will attain its least negative value at $x_i=-p^2_i$ at optimality.  Therefore, Problem (\ref{prob:reform1}) and Problem (\ref{prob:orig}) are equivalent and consequently, we can deal with Problem (\ref{prob:reform1}) from this point on.

We now apply relaxations to Problem (\ref{prob:reform1}) to address the nonlinear constraints \eqref{prob:reform1:const2} and \eqref{prob:reform1:const3} and transform the problem into a more tractable linear program.  First, we replace the quadratic term $-p_i^2$
by its piecewise linear approximation by discretizing its domain $p_i$ into $r$ segments:
\[
l_i=\bar{p}_{i0} \leq \bar{p}_{i1} \leq \ldots \leq \bar{p}_{ir}=u_i.
\]
Then at the $k$th knot:
\begin{equation*}
\begin{array}{rcl}
p_i^2 & \approx & -\bar{p}_{ik}^2 + (-2\bar{p}_{ik})(p_i - \bar{p}_{ik}) \\
& = & \bar{p}_{ik}^2 - 2\bar{p}_{ik}p_i.
\end{array}
\end{equation*}
Therefore, we can relax constraint~\eqref{prob:reform1:const2} with the following set of $r$ supporting hyperplanes of $p_i^2$:
\[
x_i \leq \bar{p}_{ik}^2 - 2 \bar{p}_{ik}p_i \quad i\in C, k\in \{1,\ldots,r\}.
\]

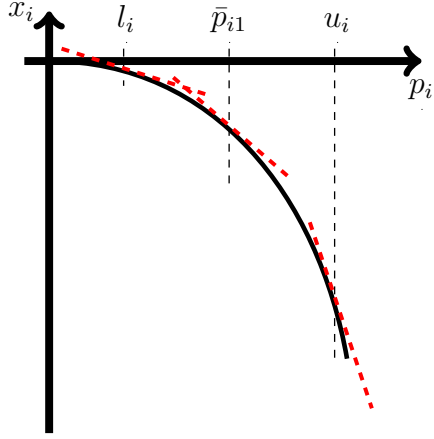
\begin{figure}[!h]
	\begin{center}

\definecolor{italian-red}{rgb}{0.965,0,0}
\definecolor{italian-green}{rgb}{0,0.45,0}

\begin{tikzpicture}[scale=0.33]

\draw[line width=3.0pt,->](-1,0)--(15,0);
\draw[line width=.1pt,smooth] (15,-1.95)--(15,-1.95) node[anchor=south] {{\bf $p_i$}};
\draw[line width=3.0pt,->](0,-15)--(0,2) node[anchor=east] {{\bf $x_i$}};

\draw[line width = 1.75pt, smooth] (0,0) to [out=0, in=100] (12, -12);

\draw[line width=.5pt,dashed] (3.,.5)--(3.,-1);
\draw[line width=.1pt,smooth] (3.1,.75)--(3.1,.75) node[anchor=south] {{$l_i$}}; 
\draw[line width=.5pt,dashed] (11.5,.5)--(11.5,-12);
\draw[line width=.1pt,smooth] (11.6,.75)--(11.6,.75) node[anchor=south] {{$u_i$}}; 




\draw[line width=1.5pt,color = red, dashed](.5,.5)--(6.5,-1.4);


\draw[line width=.5pt,dashed] (7.25,.5)--(7.25,-5);
\draw[line width=.1pt,smooth] (7.25,.75)--(7.25,.75) node[anchor=south] {{$\bar{p}_{i1}$}}; 
\draw[line width=1.5pt,color = red, dashed](5.,-.66)--(9.8,-4.8);


\draw[line width=1.5pt,color = red, dashed](10.5,-6.5)--(13,-14);


 			
%
\end{tikzpicture}
	\end{center}
	\caption{\footnotesize Linearizations of $x_i=-p_i^2$ at three points $l_i$, $\bar{p}_{i1}$ and $u_i$: dashed red
		lines are the linearizations
		\label{fig:lin}}
\end{figure}

Figure~\ref{fig:lin} illustrates how this linearization technique to approximate the quadratic function $x_i=-p_i^2$.

Next, for constraint~\eqref{prob:reform1:const3}, we relax the bilinear terms by introducing McCormick envelopes (\cite{mccormick:MP:75, Al-Khayyal:MOR:83}) as follows:
\begin{align*}
& y_{ij} \geq l_jp_i + l_ip_j - l_il_j, \\
& y_{ij} \leq u_jp_i + l_ip_j - l_iu_j, \\
& y_{ij} \leq l_jp_i + u_ip_i - l_ju_i, \\
& y_{ij} \geq u_jp_i + u_ip_j - u_iu_j.
\end{align*}

With the above relaxations, Problem~\eqref{prob:reform1} becomes:
\begin{subequations}
\label{prob:reform2}
\begin{align}
		\label{prob:reform2:obj} \max\limits_{p,x,y} \quad & \sum\limits_{i\in I}{ \left( \alpha_{i} p_i - \beta_{i} x_i +
	\sum_{j\in C_i}{\crosse_{ij} y_{ij}} \right) } \vspace{1mm} & \\
		\label{prob:reform2:const1} \st \quad & Ap \leq b, & \vspace{1mm}  \\
		\label{prob:reform2:const2} & x_i \leq \bar{p}_{ik}^2 - 2\bar{p}_{ik}p_i & i\in I, k\in \{1,\ldots,r\},  \\
		\label{prob:reform2:const31} & y_{ij} \geq l_jp_i + l_ip_j - l_il_j, & i\in I, j \in C_i  \\
		\label{prob:reform2:const32} & y_{ij} \leq u_jp_i + l_ip_j - l_iu_j, & i\in I, j \in C_i  \\
		\label{prob:reform2:const33} & y_{ij} \leq l_jp_i + u_ip_i - l_ju_i, & i\in I, j \in C_i  \\
		\label{prob:reform2:const34} & y_{ij} \geq u_jp_i + u_ip_j - u_iu_j, & i\in I, j \in C_i  \\
		\label{prob:reform2:const4} & l \leq p \leq u.
\end{align}
\end{subequations}

\section{Uncertainty}
\label{sec:uncertainty}
In our MMPO model, we assume that the uncertainty is on $\beta_{i}$ and $\gamma_{ij}$ which vary between their corresponding
lower and upper bound values. Our goal is to find optimal price decisions so that the decisions remain optimal
with respect to all possible realizations of $\beta$ and $\gamma$.
In this section, we first propose an RO approach to achieve the above goal. We then show how we construct an uncertainty set in building an RO model. Finally, we describe our method to solve the
RO model.

\subsection{The Robust Optimization Model for MPPO}
The RO model is constructed by considering the worst-case scenario under all possible realizations of the uncertain parameter. Therefore, the RO counterpart of the objective function of Problem~\eqref{prob:reform2} is written as a {\it Max-Min} problem as follows:

\begin{align*}
&& \max\limits_{p,x,y} \quad & \min_{(\beta,\gamma) \in \uset}\sum\limits_{i\in I}{ \left( \alpha_{i} p_i - \beta_{i} x_i + \sum_{j\in C_i}{\crosse_{ij} y_{ij}} \right) } \\
&=& \max\limits_{p,x,y} \quad & \left( \sum\limits_{i\in I}{ \alpha_{i} p_i } + \min_{(\beta,\gamma) \in \uset}\sum\limits_{i\in I}{ \Bigl( -\beta_{i} x_i + \sum_{j\in C_i}{\crosse_{ij} y_{ij}} \Bigr) } \right) \\
\end{align*}
Hence, Problem~\eqref{prob:reform2} under uncertainty becomes:
\begin{equation}
 \label{prob:reform3}
 \begin{array}{rl}
  \max\limits_{p,x,y} \quad & \left( \sum\limits_{i\in I}{ \alpha_{i} p_i } + \min\limits_{(\beta,\gamma) \in \uset}\sum\limits_{i\in I}{ \Bigl( -\beta_{i} x_i + \sum\limits_{j\in C_i}{\crosse_{ij} y_{ij}} \Bigr) } \right) \vspace{3mm} \\
  \st \quad & \text{\eqref{prob:reform2:const1} - \eqref{prob:reform2:const4}},
 \end{array}
\end{equation}
where $\uset$ is the uncertainty set of $(\beta, \gamma)$. We will describe a method for constructing the uncertainty set later.

By introducing an ancillary variable $\eta$, we reformulate the model~\eqref{prob:reform3} as
\begin{subequations}
 \label{prob:reform4}
 \begin{align}
  \max\limits_{p,x,y} \quad & \eta+\sum\limits_{i\in I}{ \alpha_{i} p_i  }  & \nonumber \vspace{3mm} \\
  \st \quad & \eta \leq \min\limits_{(\beta,\gamma)\in \uset}\left( \sum\limits_{i\in I}{ \Bigl( -\beta_{i} x_i + \sum_{j\in C_i}{\crosse_{ij} y_{ij}} \Bigr) } \right)  & \label{prob:reform4:const1}  \vspace{3mm} \\
   & \text{\eqref{prob:reform2:const1} - \eqref{prob:reform2:const4}}. \nonumber
 \end{align}
\end{subequations}

\subsection{Uncertainty Set Construction}
We denote $S_i^{RO}$ and $T_{ij}^{RO}$ as two pre-determined numbers of uniformly distributed realizations $\hc_{is}$ and $\hat{\gamma}_{ijs}$
for uncertain parameters
$\beta_i$ and $\gamma_{ij}$, respectively, from their corresponding interval estimates.
Given the set of realizations, we define the following attributes for the uncertainty set for
each product $i$:
\BI
	\item Mean of scenarios:
	\[
		\bar{\beta}_i = \frac{\sum\limits_{s=1}^{S_i^{RO}}{\hc_{is}}}{S_i^{RO}}.
	\]
	$\bar{\beta}_i$ can be interpreted as the nominal value for the sensitivity parameter $\beta_i$.
	Likewise, for $\gamma_{ij}$, we have:
	\[
	\bar{\gamma}_{ij} = \frac{\sum\limits_{s=1}^{T_{ij}^{RO}}{\hat{\gamma}_{is}}}{T_{ij}^{RO}}.
	\]
	
	\item Standard deviation of scenarios:
The standard deviations of the samples of realizations for $\beta$ and $\gamma$ are calculated as follows.
	\[
		\tilde{\beta}_i = \sqrt{\frac{\sum\limits_{s=1}^{S_i^{RO}}{(\hc_{is} - \bar{\beta}_i)^2}}{S_i^{RO}-1}},
	\]
	and
	\[
		\tilde{\gamma}_{ij} = \sqrt{\frac{\sum\limits_{s=1}^{T_{ij}^{RO}}{(\hat{\gamma}_{ijs} - \bar{\gamma}_{ij})^2}}{T_{ij}^{RO}-1}}.
	\]
\EI

Now, for $k = \sum_{i\in I}{|C_i|}$ and writing $\boldsymbol\beta$ for the vectors of all $\beta$ and $\boldsymbol\gamma \in \Real^k$ for the vector of all $\gamma_{ij}$,
we construct the uncertainty set $\uset_{\Delta}$ as follows:

\begin{subequations}
\begin{align}
\uset_{\Delta} = \Bigl\{
\left[
\begin{array}{c}
\boldsymbol\beta \\
\boldsymbol\gamma
\end{array}
\right]  \in & \Real^{n+k} :\Bigr. \nonumber \\
	&  \max\limits_{i\in I}{\frac{|\beta_i - \bar{\beta}_i|}{\tilde{\beta}_i}} \leq \Delta, \label{unset:1}\\
	& \hc^{min}_{i} \leq \beta_i \leq \hc^{max}_{i}, \ \forall i\in I, \label{unset:2}\\
	& \max\limits_{i\in I, j\in C_i}{\frac{|\gamma_{ij} - \bar{\gamma}_{ij}|}{\tilde{\gamma}_{ij}}} \leq \Delta, \label{unset:3}\\
	& \hat{\gamma}^{min}_{ij} \leq \gamma_{ij} \leq \hat{\gamma}^{max}_{ij}, \ \forall i\in I,j\in C_i \Bigl. \Bigr\}, \label{unset:4}
\end{align}
\end{subequations}
where $\hc^{min}_{i} = \min_{s \in \{1,\ldots,S_i^{RO}\}}\{\hc_{is}\}$, $\hc^{max}_{i} = \max_{s \in \{1,\ldots,S_i^{RO}\}}\{\hc_{is}\}$,
$\hat{\gamma}^{min}_{ij} = \min_{s \in \{1,\ldots,T_{ij}^{RO}\}}\{\hat{\gamma}_{ijs}\}$ and
$\hat{\gamma}^{max}_{ij} = \max_{s \in \{1,\ldots,T_{ij}^{RO}\}}\{\hat{\gamma}_{ijs}\}$.

Note that the constraints~\eqref{unset:1} and \eqref{unset:3} in the above construction ensure that the maximum normalized deviation of the uncertain parameter $\beta_i$ and $\gamma_{ij}$, respectively, from the
nominal values among all the products remains under a certain positive scalar $\Delta$.
The parameter $\Delta$ is called {\it budget of uncertainty} \cite{bertsimas:TOR:06}. Constraintd
\eqref{unset:2} and \eqref{unset:4} also define bounds on individual $\beta_i$ and $\gamma_{ij}$, respectively.

The budget of uncertainty constraint could  also be specified by imposing a limit on the {\it total} amount of deviation across all
products \cite{bertsimas:ieee:13} as follows:
\[
\sum\limits_{i\in I}{\frac{|\beta_i - \bar{\beta}_i|}{\tilde{\beta}_i}} +
\sum_{i\in I}{\sum_{j\in C_i}{\frac{|\gamma_{ij} - \bar{\gamma}_{ij}|}{\tilde{\gamma}_{ij}}}}\leq \Delta.
\]
We therefore reformulate the uncertainty set as a polytope:
\begin{equation}
\label{set:uncertset1}
\begin{array}{rll}
\ruset_{\Delta} = \Bigl\{
\left[
\begin{array}{c}
\boldsymbol\beta \\
\boldsymbol\gamma \\
\boldsymbol\sigma^\beta \\
\boldsymbol\sigma^\gamma
\end{array}
\right]
& \hspace{-3mm} \in \Real^{2n+2k} :\Bigr.  & \\
& \hspace{-5mm} \sigma_i^{\beta} \leq \Delta, & \forall i\in I, \vspace{2mm} \\
& \hspace{-5mm} \sigma_{ij}^{\gamma} \leq \Delta, & \forall i\in I, j\in C_i, \vspace{2mm} \\
& \hspace{-5mm} - \tilde{\beta}_i \cdot \sigma_i^{\beta} \leq \beta_i - \bar{\beta}_i \leq \tilde{\beta}_i
\cdot \sigma_i^{\beta}, & \forall i\in I, \vspace{2mm} \\
& \hspace{-5mm} \hc^{min}_{i} \leq \beta_i \leq \hc^{max}_{i}, & \forall i\in I, \\
& \hspace{-5mm} - \tilde{\gamma}_{ij} \cdot \sigma_{ij}^{\gamma} \leq \gamma_{ij} - \bar{\gamma}_{ij} \leq \tilde{\gamma}_{ij} \cdot \sigma_{ij}^{\gamma}, & \forall i\in I, j\in C_i, \vspace{2mm} \\
& \hspace{-5mm} \hat{\gamma}^{min}_{i} \leq \gamma_{ij} \leq \hat{\gamma}^{max}_{ij}, & \forall i\in I, j\in C_i   \Bigl. \Bigr\}
\end{array}
\end{equation}

\subsection{Solution Approach for the RO Model}
Note that the right-hand side of the constraint~\eqref{prob:reform4:const1} constitutes an optimization problem. We substitute the uncertainty set $\uset$ in ~\eqref{prob:reform4:const1} with a hyperplane derived from \eqref{set:uncertset1}, and yield the following optimization problem:

\begin{subequations}
\label{prob:ro}
\begin{align}
\min\limits_{\beta,\gamma, \sigma^{\beta}, \sigma^{\gamma}} \quad & \sum\limits_{i\in I}{ \Bigl( -\beta_{i} x_i + \sum_{j\in C_i}{\crosse_{ij} y_{ij}} \Bigr)} & \nonumber \\
\st \quad & \sigma_i^{\beta} \leq \Delta, & \forall i\in I, \label{prob:ro:constr1} \\
& \sigma_{ij}^{\gamma} \leq \Delta, & \forall i\in I, j\in C_i, \label{prob:ro:constr2}\\
&  \beta_i - \bar{\beta}_i \leq \tilde{\beta}_i
\cdot \sigma_i^{\beta}, & \forall i\in I, \vspace{2mm} \label{prob:ro:constr3}\\
& - \tilde{\beta}_i \cdot \sigma_i^{\beta} \leq \beta_i - \bar{\beta}_i,
& \forall i\in I, \vspace{2mm} \label{prob:ro:constr4 }\\
& \hc^{min}_{i} \leq \beta_i  & \forall i\in I, \label{prob:ro:constr5}\\
& \beta_i \leq \hc^{max}_{i},  & \forall i\in I, \label{prob:ro:constr6}\\
&  \gamma_{ij} - \bar{\gamma}_{ij} \leq \tilde{\gamma}_{ij} \cdot \sigma_{ij}^{\gamma}, & \forall i\in I, j\in C_i, \vspace{2mm} \label{prob:ro:constr7}\\
& - \tilde{\gamma}_{ij} \cdot \sigma_{ij}^{\gamma} \leq \gamma_{ij} - \bar{\gamma}_{ij}, & \forall i\in I, j\in C_i, \vspace{2mm} \label{prob:ro:constr8}\\
& \hat{\gamma}^{min}_{i} \leq \gamma_{ij}  & \forall i\in I, j\in C_i \label{prob:ro:constr9}\\
& \gamma_{ij} \leq \hat{\gamma}^{max}_{ij}, & \forall i\in I, j\in C_i \label{prob:ro:constr10} \\
& \sigma_i^{\beta}, \sigma_{ij}^{\gamma} \geq 0, & i\in I, j\in C_i. \nonumber
\end{align}
\end{subequations}

By introducing dual multipliers $\mu^{\beta}_i$, $\mu^{\gamma}_{ij}$,
$\pi_{1i}^{\beta}$, $\pi_{2i}^{\beta}$,
$\lambda_{1i}^{\beta}$, $\lambda_{2i}^{\beta}$,
$\pi_{1ij}^{\gamma}$, $\pi_{2ij}^{\gamma}$,
$\lambda_{1ij}^{\gamma}$ and $\lambda_{2ij}^{\gamma}$
associated with constraints \eqref{prob:ro:constr1}-\eqref{prob:ro:constr10}, respectively, the dual problem becomes:
\begin{equation}
\label{prob:dualreform}
\begin{array}{rll}
\hspace{-5mm} \max
 & \sum\limits_{i\in I}{\bigl( -\Delta \mu_i^{\beta} - \bar{\beta}_{i} \pi_{1i}^{\beta}
	+ \bar{\beta}_i \pi_{2i}^{\beta} + \hc^{min}_{i} \lambda_{1i}^{\beta}
	-  \hc^{max}_{i} \lambda_{2i}^{\beta}  \bigr)} & \\
& + \sum\limits_{i\in I}{\sum\limits_{j\in C_i}{\Bigl( -\Delta \mu_{ij}^{\gamma} - \bar{\gamma}_{i} \pi_{1ij}^{\gamma}
	+ \bar{\gamma}_i \pi_{2ij}^{\gamma} + \hat{\gamma}^{min}_{ij} \lambda_{1ij}^{\gamma} \Bigl.}} \\
& \hspace{19mm}-  \hat{\gamma}^{max}_{ij} \lambda_{2ij}^{\gamma} \Bigr. \Bigl) & \vspace{2mm} \\
\st & -\pi_{1i}^{\beta} + \pi_{2i}^{\beta} + \lambda_{1i}^{\beta} - \lambda_{2i}^{\beta} = -x_i & i \in I,\vspace{2mm} \\
& -\mu_{i}^{\beta} + \tilde{\beta}_{i} \pi_{1i}^{\beta} + \tilde{\beta}_{i} \pi_{2i}^{\beta} \leq 0, & i \in I, \vspace{2mm} \\
& -\pi_{1ij}^{\gamma} + \pi_{2ij}^{\gamma} + \lambda_{1ij}^{\gamma} - \lambda_{2ij}^{\gamma} = y_{ij} & i \in I, j\in C_i, \vspace{2mm} \\
& -\mu_{ij}^{\gamma} + \tilde{\gamma}_{ij} \pi_{1ij}^{\gamma} + \tilde{\gamma}_{ij} \pi_{2ij}^{\gamma} \leq 0, & i \in I, j\in C_i, \vspace{2mm} \\
& \mu_i^{\beta}, \pi_{1i}^{\beta}, \pi_{2i}^{\beta}, \lambda_{1i}^{\beta}, \lambda_{2i}^{\beta} \geq 0. & i \in I, \vspace{2mm} \\
& \mu_{ij}^{\gamma}, \pi_{1ij}^{\gamma}, \pi_{2ij}^{\gamma}, \lambda_{1ij}^{\gamma}, \lambda_{2ij}^{\gamma} \geq 0. & i \in I, j\in C_i.
\end{array}
\end{equation}

In accordance with weak duality, any feasible solution to \eqref{prob:dualreform} provides a lower bound
to \eqref{set:uncertset1}. Therefore, we can enforce  \eqref{prob:reform4:const1} by adding the following  constraint
\begin{align*}
\eta \leq & \sum\limits_{i\in I}{\bigl( -\Delta \mu_i^{\beta} - \bar{\beta}_{i} \pi_{1i}^{\beta}
	+ \bar{\beta}_i \pi_{2i}^{\beta} + \hc^{min}_{i} \lambda_{1i}^{\beta}
	-  \hc^{max}_{i} \lambda_{2i}^{\beta}  \bigr) + } \\
& \sum\limits_{i\in I}{\sum\limits_{j\in C_i}{\Bigl( -\Delta \mu_{ij}^{\gamma} - \bar{\gamma}_{i} \pi_{1ij}^{\gamma}
		+ \bar{\gamma}_i \pi_{2ij}^{\gamma} + \hc^{min}_{i} \lambda_{1ij}^{\gamma} - \hc^{max}_{i} \lambda_{2ij}^{\gamma} \Bigr)}}
\end{align*}
where $(\mu^{\beta}, \pi_{1\cdot}^{\beta}, \pi_{2\cdot}^{\beta}, \lambda_{1\cdot}^{\beta},
\lambda_{2\cdot}^{\beta}, \mu^{\gamma}, \pi_{1}^{\gamma}, \pi_{2}^{\gamma}, \lambda_{1}^{\gamma},\lambda_{2}^{\gamma})$ is feasible solution for \eqref{prob:dualreform}. Therefore,
the {\it Max-Min } problem can be reformulated as:
\begin{subequations}
	\small
	\label{prob:reform5}
	\begin{align}
	\hspace{-5mm} Z_{\Delta}^*=\max & \hspace{1mm} \eta + \sum\limits_{i\in I}{ \alpha_{i} p_i } & \vspace{3mm} \\
	\st & \hspace{1mm} \eta \leq  \sum\limits_{i\in I}{\bigl( -\Delta \mu_i^{\beta} - \bar{\beta}_{i} \pi_{1i}^{\beta}
		+ \bar{\beta}_i \pi_{2i}^{\beta} + \hc^{min}_{i} \lambda_{1i}^{\beta}
		 \bigr. } &  \nonumber \vspace{3mm} \\
	& \hspace{5mm} -  \hc^{max}_{i} \lambda_{2i}^{\beta} \bigl. \bigr) +\sum\limits_{i\in I}{\sum\limits_{j\in C_i}{\Bigl( -\Delta \mu_{ij}^{\gamma} - \bar{\gamma}_{i} \pi_{1ij}^{\gamma}
			+ \bar{\gamma}_i \pi_{2ij}^{\gamma} \Bigr.}} & \nonumber \\
	& \hspace{5mm} \Bigl. + \hat{\gamma}^{min}_{ij} \lambda_{1ij}^{\gamma} - \hat{\gamma}^{max}_{ij} \lambda_{2ij}^{\gamma} \Bigr) & \\
	& -\pi_{1i}^{\beta} + \pi_{2i}^{\beta} + \lambda_{1i}^{\beta} - \lambda_{2i}^{\beta} = -x_i & i \in I,\vspace{2mm} \\
	& -\mu_{i}^{\beta} + \tilde{\beta}_{i} \pi_{1i}^{\beta} + \tilde{\beta}_{i} \pi_{2i}^{\beta} \leq 0, & i \in I, \vspace{2mm} \\
	& -\pi_{1ij}^{\gamma} + \pi_{2ij}^{\gamma} + \lambda_{1ij}^{\gamma} - \lambda_{2ij}^{\gamma} = y_{ij} & i \in I, j\in C_i, \vspace{2mm} \\
	& -\mu_{ij}^{\gamma} + \tilde{\gamma}_{ij} \pi_{1ij}^{\gamma} + \tilde{\gamma}_{ij} \pi_{2ij}^{\gamma} \leq 0, & i \in I, j\in C_i, \vspace{2mm} \\
	& \text{\eqref{prob:reform2:const1} - \eqref{prob:reform2:const4}}, & \\
	& \mu_i^{\beta}, \pi_{1i}^{\beta}, \pi_{2i}^{\beta}, \lambda_{1i}^{\beta}, \lambda_{2i}^{\beta} \geq 0. & i \in I, \vspace{2mm} \\
	& \mu_{ij}^{\gamma}, \pi_{1ij}^{\gamma}, \pi_{2ij}^{\gamma}, \lambda_{1ij}^{\gamma}, \lambda_{2ij}^{\gamma} \geq 0. & i \in I, j\in C_i.		
	\end{align}
\end{subequations}

The above formulation constitutes the RO model associated with a budget of uncertainty parameter $\Delta$. This RO model is a conventional linear programming model that can be readily solved using standard linear programming solvers. We report the results of our numerical experiments based on this model in the subsequent section.

\section{Numerical Experiments}
\label{sec:NumExam}
We implemented the RO model presented above using Python 3.5.2. We employed Gurobi 7.0 \cite{gurobi:7}
to solve the linear programming model in \eqref{prob:reform5}. All tests reported were conducted on Amazon's
Elastic Compute Cloud running Amazon Linux AMI 2018.03 on a machine with dual Intel(R) Xeon(R)
16-core E5-2686 v4 CPU @ 2.3GHz and 500GB memory.

\subsection{Test Instance Description}
We tested the proposed approach on 40 different instances from the data of a major food distributor
in the United States. The company sells products to its customers via its geographically scattered
divisions. Table~\ref{tbl:inststats} lists
the name of the test instances (such that each instance represents a specific division), the number of products and number of cross terms (nonzero $\gamma_{ij}$ terms) within a given time period in each instance.

\FloatBarrier
\begin{table}[!ht]
	\footnotesize
	\begin{center}
\begin{tabular} {lrr|lrr}
	\hline\hline
	Sample & \# Products & \# Cross Terms & Sample & \# Products & \# Cross Terms \\
	\hline
		S1 & 14,072 & 153 & S2 & 11,563 & 247 \\
		S3 & 15,665 & 151 & S4 & 1,148 & 19 \\
		S5 & 310 & 5 & S6 & 1,975 & 23 \\
		S7 & 887 & 16 & S8 & 987 & 15 \\
		S9 & 3,934 & 23 & S10 & 2,078 & 35\\
		S11 & 1,635 & 14 & S12 & 2,008 & 16 \\
		S13 & 1,281 & 18 & S14 & 629 & 10\\
		S15 & 11,442 & 149 & S16 & 2,940 & 79 \\
		S17 & 1,807 & 70 & S18 & 924 & 53 \\
		S19 & 4,222 & 126 & S20 & 1,675 & 26 \\
		S21 & 1,768 & 11 & S22 & 333 & 22 \\
		S23 & 3,353 & 29 & S24 & 1,232 & 34 \\
		S25 & 389,954 & 0 & S26 & 288,719 & 0 \\
		S27 & 69,107 & 1,313 & S28 & 141,475 & 3,409 \\
		S29 & 35,679 & 0 & S30 & 298,750 & 0\\
		S31 & 23,269 & 439 & S32 & 9,518 & 201 \\
		S33 & 23,092 & 413 & S34 & 23,249 & 481 \\
		S35 & 20,593 & 366 & S36 & 18,862 & 232 \\
		S37 & 375,069 & 0 & S38 & 271,054 & 0\\
		S39 & 83,211 & 2,139 & S40 & 112,555 & 2,686 \\
			\hline\hline
\end{tabular} 	
\end{center}
\caption{\footnotesize Statistics about instances: Sample name and number of products
	and number of cross terms \label{tbl:inststats} }
\end{table}

\subsection{Parameter Values of the RO method}
\label{subsec:valuedelta}
The value of $\Delta$ plays a crucial role in our tests. We select and test the values as follows.
First, let us revisit
the constraint involving $\Delta$:
\begin{align*}
& \max_{i\in I}\left\{\frac{|\beta_i - \bar{\beta}_i|}{\tilde{\beta}_i}\right\} \leq \Delta, \quad \forall i\in I,\\
& \max\limits_{i\in I, j\in C_i}{\frac{|\gamma_{ij} - \bar{\gamma}_{ij}|}{\tilde{\gamma}_{ij}}} \leq \Delta, \quad \forall i \in I, j\in C_i.
\end{align*}
Obviously, the minimum value for $\Delta$ is zero, in which case all the $\beta_i$ values will be equal to their
nominal value $\bar{\beta}_i$. We can also find the maximum value for $\Delta$ by:
\[
\Delta_{max}^{\beta} = \max_{i\in I}\left\{ \max\left\{\frac{|\hc^{min}_{i} - \bar{\beta}_i|}{\tilde{\beta}_i},
\frac{|\hc^{max}_{i} - \bar{\beta}_i|}{\tilde{\beta}_i} \right\} \right\},
\]
and
\[
\Delta_{max}^{\gamma} = \max_{i\in I, j\in C_i}\left\{ \max\left\{\frac{|\hat{\gamma}^{min}_{i} - \bar{\gamma}_i|}{\tilde{\gamma}_i},
\frac{|\hat{\gamma}^{max}_{i} - \bar{\gamma}_i|}{\tilde{\gamma}_i} \right\} \right\}.
\]
Therefore, we choose $\Delta_{max} = \max\{\Delta_{max}^{\beta}, \Delta_{max}^{\gamma}\}$.
We then divide the interval [0,$\Delta_{max}$] into 10 equally distanced segments, namely $[\Delta_1, \Delta_2]$, $[\Delta_2, \Delta_3], \ldots, [\Delta_{10}, \Delta_{11}]$ where
$\Delta_1=0$ and $\Delta_{11}=\Delta_{max}$, and test each of the $\Delta$ values in our subsequent experiments.

Other parameters of the model are $S_i^{RO}$ and $T_{ij}^{RO}$, the numbers of samples used for calculating the means and variances for $\beta_i$ and $\gamma_{ij}$.
We set these two parameters to 50 and 5, respectively, reflecting the fact that the cross price elasticities are usually associated with smaller intervals compared with the own price elasticities, as found in our real data sets.

\subsection{Test Results}
To present our results, we show how changing the value of $\Delta$ affects the optimal value of the robust counterpart
problem~\eqref{prob:reform5}. To do this, we calculate the percent difference between the maximum and the minimum value of $Z^*$
associated with values of $\Delta$ in problem~\eqref{prob:reform5}. The maximum and minimum values of $Z^*$ are
attained when $\Delta = 0$ and $\Delta = \Delta_{max}$, i.e., the least conservative case and the most
conservative case, respectively. We summarize the results in Table~\ref{tbl:declines}
by only displaying percent of decline from the least conservative to the most conservative case.
We observe that the differences are all within in the range $(0.05\%, 6.96\%)$ with an average of 1.14\%.

\begin{table}[!h]
	\begin{center}
		\small
		\begin{tabular}{lr|lr}
			\hline\hline
Sample & Decline (\%) & Sample & Decline (\%) \\
\hline
S1 & 4.0869 & S2 & 0.7784 \\
S3 & 6.1058 & S4 & 5.2174 \\
S5 & 0.0964 & S6 & 0.2366 \\
S7 & 0.0844 & S8 & 0.1391 \\
S9 & 3.5418 & S10 & 0.1005 \\
S11 & 0.3744 & S12 & 0.1979 \\
S13 & 0.1200 & S14 & 0.1398 \\
S15 & 0.3392 & S16 & 0.0951 \\
S17 & 0.1019 & S18 & 0.0716 \\
S19 & 0.1092 & S20 & 0.1906 \\
S21 & 0.0898 & S22 & 0.0757 \\
S23 & 0.2508 & S24 & 0.0681 \\
S25 & 0.8945 & S26 & 1.8451 \\
S27 & 0.0501 & S28 & 0.0691 \\
S29 & 2.4146 & S30 & 0.5939 \\
S31 & 6.9645 & S32 & 2.6059 \\
S33 & 1.2579 & S34 & 0.2168 \\
S35 & 0.1478 & S36 & 0.8719 \\
S37 & 3.0746 & S38 & 1.1271 \\
S39 & 0.4145 & S40 & 0.3026 \\
\hline \hline
\end{tabular}
\caption{\footnotesize Percent decline from the least conservative to the most conservative solutions with RO approach. \label{tbl:declines}}
\end{center}
\end{table}

The distribution of the percent differences is depicted in the box plot of Figure~\ref{boxplot}.
The detailed results are shown in the Appendix (Table~\ref{tbl:deltaoptval}).

\FloatBarrier
\begin{figure}
	\begin{center}
		\includegraphics[width=.75\linewidth]{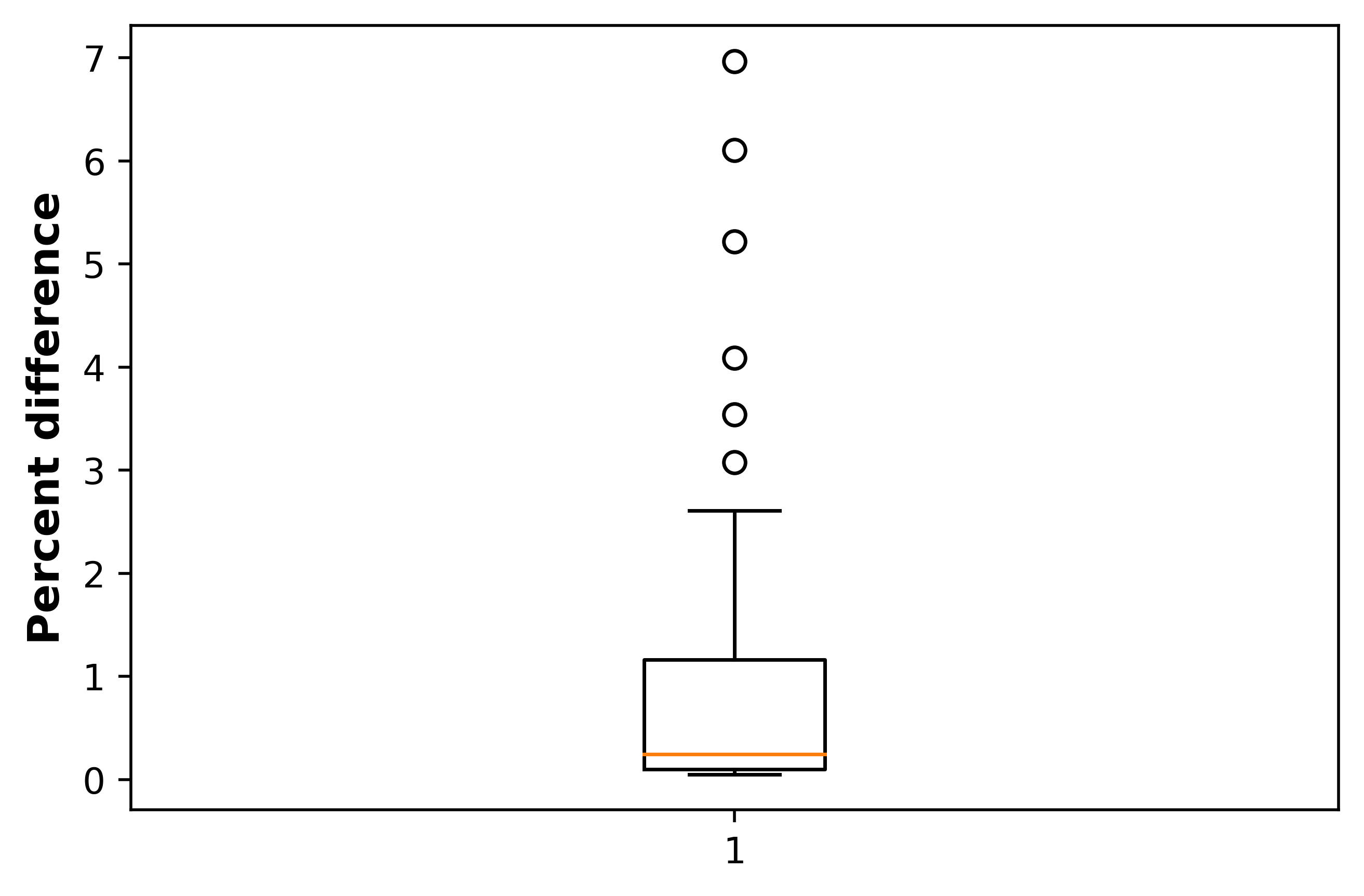}
	\end{center}
\vspace{-5mm}
\caption{\footnotesize  Boxplot of percent difference between least and most conservative solutions of RO problem for instances.\label{boxplot}}
\end{figure}

\subsection{Simulation Validation}
We further validate our RO approach by comparing results from a set of Monte Carlo simulation tests with our RO solutions.
We run a Monte Carlo simulations by taking $S^{sim}$ uniformly distributed samples $\beta_{is}$ from
the interval $[\hc^{min}_{i}, \hc^{max}_{i}]$ and $\gamma_{ijs}$ from
$[\hat{\gamma}_{ij}^{min},\hat{\gamma}_{ij}^{max}]$,
and then solving each of the price optimization problems with given sampling values for $\beta$ and $\gamma$:
\begin{equation}
\label{prob:sim}
\begin{array}{rl}
f(\boldsymbol\beta_{s},\boldsymbol\gamma_s) = \max\limits_{p,x,y} & \sum\limits_{i\in I}{ ( \alpha_{i} p_i - \beta_{is} x_i +
	\sum_{j\in C_i}{\crosse_{ij} y_{ij}} ) } \vspace{1mm} \\
		\st & \eqref{prob:reform2:const1}-\eqref{prob:reform2:const4}.
\end{array}
\end{equation}
The value of the simulation will then be the mean of optimal values from each replication, i.e. $Z_{sim}^* = \frac{1}{S^{sim}} \sum\limits_s{f(\boldsymbol\beta_{s},\boldsymbol\gamma_s)}$.

In these experiments, with $S^{sim}=50$, we compare the values of the simulation $Z^*_{sim}$ with the values of RO formulation \eqref{prob:reform5}.
This comparison discloses the degree of conservativeness of the RO solutions. To perform such a comparison, we compare
$Z^*_{sim}$ with the least conservative case (achieved at $\Delta_1$), the most conservative case $\Delta_{max}$
and the average across all $\Delta$ values, and calculate the corresponding difference ratios. For instance the difference ratio in the column under $\Delta_1$ is attained from
\[
	\frac{Z^*_{sim} - Z^*_{\Delta_1}}{Z^*_{\Delta_1}}.
\]
Table~\ref{tbl:declines1} summarizes these comparisons. while a negative value of difference ratio in the table implies that the RO value is better
than the simulation,  a positive value shows conversely. Note that most of the instances (35 out of 40) have positive values in $\Delta_{max}$ column, indicating RO results are
worse than the simulation values. However, in most of such cases (30 out of 35) the difference ratios are less than 0.5\%, and there are only three cases where the difference ratios are over 1\%. This implies that for
the most conservative scenarios, RO solution yields results with little losses as compared to the solutions from the simulation.
On the other hand, in the least conservative case (under column $\Delta_1$), most (31 of 40) differences are negative, demonstrating
RO results are better than those from the simulations, while there are only nine instances where it is not and only in one of which the difference ratio is above 0.32\%.

\FloatBarrier
\begin{table}[!h]
	\begin{center}
		\scriptsize
		\begin{tabular}{lccc|lccc}
		\hline \hline
Sample & $\Delta_1$ (\%) & Average (\%) & $\Delta_{max}$ (\%) & Sample & $\Delta_1$ (\%) & Average (\%) & $\Delta_{max}$ (\%) \\
			\hline
S1 & -1.8890 & 0.2668 & 0.9437 & S2 & -0.0190 & 0.2335 & 0.2986 \\
S3 & -2.8905 & 0.4543 & 1.9151 & S4 & 0.2826 & 3.6168 & 5.0061 \\
S5 & -0.1005 & -0.1198 & -0.0184 & S6 & -0.1459 & 0.0224 & 0.0721 \\
S7 & 0.0542 & 0.0965 & 0.1128 & S8 & 0.1286 & 0.1395 & 0.1465 \\
S9 & -5.0707 & -2.6889 & -1.8811 & S10 & 0.0430 & 0.1518 & 0.1716 \\
S11 & -0.2461 & -0.0556 & 0.0774 & S12 & -0.0706 & 0.0239 & 0.1493 \\
S13 & -0.0518 & 0.1993 & 0.2539 & S14 & -0.0023 & 0.0169 & 0.0832 \\
S15 & -0.1280 & 0.0184 & 0.2112 & S16 & -0.0800 & 0.0541 & 0.1453 \\
S17 & -0.1146 & -0.0004 & 0.0926 & S18 & 0.0192 & 0.0281 & 0.0357 \\
S19 & -0.0865 & -0.0006 & 0.0103 & S20 & -0.1491 & -0.0081 & 0.0242 \\
S21 & -0.0953 & 0.0578 & 0.0941 & S22 & 0.0062 & 0.0327 & 0.0401 \\
S23 & -0.2302 & -0.0727 & 0.0346 & S24 & 0.0034 & 0.0293 & 0.0658 \\
S25 & -0.6957 & -0.1953 & 0.1140 & S26 & -1.5647 & -0.1403 & 0.3698 \\
S27 & -0.0236 & 0.0398 & 0.0529 & S28 & -0.1683 & -0.0924 & -0.0981 \\
S29 & -1.8245 & -03922 & 0.3490 & S30 & -0.5110 & 0.0049 & 0.1277 \\
S31 & -10.1855 & -5.0828 & -2.4091 & S32 & 1.8656 & 2.6629 & 3.4037 \\
S33 & -0.2890 & 0.2004 & 0.3969 & S34 & -0.0104 & 0.0221 & 0.0562 \\
S35 & -0.0402 & 0.0493 & 0.0873 & S36 & -1.3137 & -0.1453 & 0.296 \\
S37 & -4.7585 & -3.0283 & -2.0466 & S38 & -0.8910 & -0.236 & 0.2993 \\
S39 & -0.0165 & 0.0100 & 0.1477 & S40 & 0.3194 & 0.5114 & 0.6030 \\
			\hline \hline
		\end{tabular}
		\caption{\footnotesize Decline of the simulation optimal value
			from the RO models in the least conservative ($\Delta_1$) and the most conservative
			($\Delta_{max}$). Also decline of the simulation optimal value
			from the average of the optimal values of the RO models across all
			$\Delta$ values.\label{tbl:declines1}}
	\end{center}
\end{table}

The other major benefit that we observe in our experiments is the gain from the computational time required for our RO approach.
Table~\ref{tbl:solutiontime} lists the solution time of both RO approach and simulation. The solution times are
significantly higher in simulation comparing to those using the RO approach. This is aligned with
the expectations as the simulation needs to solve $S^{sim}$ number of linear programs whereas the RO approach needs to only solve one
linear program. The comparison of the solution time shows that the relative difference calculated by
$\frac{t_{\text{Simulation}} - t_{\text{RO}}}{t_{\text{RO}}}$ has a geometric mean of 21.2 which means
the simulation on average takes 21 times longer to solve. The advantage in computation times of our
RO approach represents a significant benefit that could be achieved in many practical applications.

\FloatBarrier
\begin{table}[!ht]
	\footnotesize
	\begin{center}
\begin{tabular} {lrr|lrr}
	\hline\hline
	Sample & RO & Simulation & Sample & RO & Simulation \\
			& \multicolumn{1}{c}{(s)} & \multicolumn{1}{c|}{(s)} & &	\multicolumn{1}{c}{(s)} 	& \multicolumn{1}{c}{(s)} \\
	\hline
S1 & 18.1 & 298.7 & S2 & 17.5 & 323.5 \\
S3 & 9.7 & 203.3 & S4 & 0.2 & 6.3 \\
S5 & 0.2 & 8.6 & S6 & 0.4 & 14.6 \\
S7 & 0.3 & 8.4 & S8 & 0.3 & 10.4 \\
S9 & 1.0 & 31.5 & S10 & 0.7 & 22.8 \\
S11 & 0.3 & 10.1 & S12 & 0.3 & 11.6 \\
S13 & 0.5 & 15.3 & S14 & 0.1 & 3.5 \\
S15 & 5.5 & 118.5 & S16 & 1.4 & 39.7 \\
S17 & 0.4 & 11.5 & S18 & 0.2 & 5.5 \\
S19 & 1.2 & 40.2 & S20 & 0.4 & 11.3 \\
S21 & 0.7 & 22.6 & S22 & 0.1 & 1.9 \\
S23 & 1.5 & 37.2 & S24 & 0.2 & 7.7 \\
S25 & 323.6 & 2746.9 & S26 & 206.0 & 1759.1 \\
S27 & 332.7 & 5622.9 & S28 & 888.8 & 35712.6 \\
S29 & 8.1 & 182.3 & S30 & 187.0 & 1842.6 \\
S31 & 106.0 & 1410.4 & S32 & 23.2 & 395.4 \\
S33 & 90.6 & 1152.9 & S34 & 79.8 & 1259.7 \\
S35 & 71.2 & 1044.9 & S36 & 29.2 & 492.9 \\
S37 & 313.9 & 2630.8 & S38 & 177.8 & 1649.9 \\
S39 & 436.6 & 12469.6 & S40 & 436.6 & 12469.6 \\
			\hline\hline
\end{tabular}
\caption{\footnotesize Solution time for RO approach and simulation in seconds \label{tbl:solutiontime}} 	
\end{center}
\end{table}

\section{Conclusions and Future Directions}
\label{sec:conclusion}

In this paper, we investigate a multiple product price optimization problem with a linear
demand function associated with prices on reference products as well as associated complementary/substitute products.
The problem is associated with uncertainty where the sensitivity parameters on the product prices are expressed
using uncertainty intervals. In contrast to the common practice of solving the problem with a
point estimator of the interval, we devise a new robust optimization based approach to handle such uncertainty intervals.

In our proposed approach, we construct an uncertainty set based on all realizations of the
parameters. The set depends on the budget of uncertainty parameter that controls the level of conservatism.
With this set, a robust optimization formulation counterpart is developed that takes the uncertainty set into account, and it can be solved as a linear programming model.
We further test our proposed model and solution approach based on several real data sets. As is expected, by increasing the budget of uncertainty value of the
uncertainty set parameter, the optimal revenue decreases. An additional set of simulation experiments shows that although
robust optimization is typically a more conservative approach, the losses of optimal values with the approach, in most instances, are quite insignificant. From the solution time perspective, while the simulation requires much longer time to solve,
the proposed RO method, designed to address the uncertainty issue, can yield robust solutions in much shorter time at the price of being
a little conservative.

We plan to further our research in several areas in the future. First, we currently use a uniform sampling approach from the intervals
of the sensitivity parameters. However, in reality, some more refined information including posterior distribution may
be collected and estimated for such parameters. A more appropriate sampling approach can be employed to better suit
the posterior distribution of the parameters with uncertainty. 	

Also, in our robust optimization approach, we formulate the uncertainty set as explained in
$\uset_{\Delta}$, which was linearized in \eqref{set:uncertset1}. However, there are some
other approaches in the literature that might result in different approximation  of the uncertainty
associated with the parameters (for example \cite{Bertsimas:MP:2018}.) Such new approximation methods are worth exploring in the future.

Finally, we find a simple linear demand response function appropriate for our application. More complicated demand response models, such
as MNL demand functions, may be more appropriate and accurate in other price optimization applications. Hence, studying RO approach with such models
could be an interesting future direction.

\bibliographystyle{hplain}
\bibliography{references}

\begin{thebibliography}{10}

\bibitem{Agrawal:MSOM:00}
V.~Agrawal and S.~Seshadri.
\newblock Impact of uncertainty and risk aversion on price and order quantity
  in the newsvendor problem.
\newblock {\em Manufacturing \& Service Operations Management}, 2:410--423, 10
  2000.

\bibitem{Al-Khayyal:MOR:83}
F.A. Al-Khayyal and J.E. Falk.
\newblock Jointly constrained biconvex programming.
\newblock {\em Mathematics of Operations Research}, 8(2):159--317, 1983.

\bibitem{Ben-Tal:1998:MOR}
A.~Ben-Tal and A.~Nemirovski.
\newblock Robust convex optimization.
\newblock {\em Mathematics of Operations Research}, 23(4):769--805, 1998.

\bibitem{Ben-tal:99:ORL}
A.~Ben-Tal and A.~Nemirovski.
\newblock Robust solutions of uncertain linear programs.
\newblock {\em Operations Research Letters}, 25:1--13, 1999.

\bibitem{Ben-Tal:00:MP}
A.~Ben-Tal and A.~Nemirovski.
\newblock Robust solutions of linear programming problems contaminated with
  uncertain data.
\newblock {\em Mathematical Programming}, 88(3):411--424, 2000.

\bibitem{Bertsimas:OR:2009}
D.~Bertsimas and D.~Brown.
\newblock Constructing uncertainty sets for robust linear optimization.
\newblock {\em Operations Research}, 57(6):1483--1495, 2009.

\bibitem{Bertsimas:SiamRev:11}
D.~Bertsimas, D.B. Brown, and C.~Caramanis.
\newblock Theory and applications of robust optimization.
\newblock {\em SIAM Review}, 53(3):464--501, 2011.

\bibitem{Bertsimas:MP:2018}
D.~Bertsimas, V.~Gupta, and N.~Kallus.
\newblock Data-driven robust optimization.
\newblock {\em Mathematical Programming}, 167(2):235--292, 2018.

\bibitem{Bertsimas:ORL:04}
D.~Bertsimas, D.~Pachamanova, and M.~Sim.
\newblock Robust linear optimization under general norms.
\newblock {\em Operations Research Letters}, 32(6):510--516, 2004.

\bibitem{Bertsimas:mp:03}
D.~Bertsimas and M.~Sim.
\newblock Robust discrete optimization and network flows.
\newblock {\em Mathematical Programming}, 98(1):49--71, 2003.

\bibitem{Bertsimas:OR:04}
D.~Bertsimas and M.~Sim.
\newblock The price of robustness.
\newblock {\em Operations Research}, 52(1):35--53, 2004.

\bibitem{bertsimas:TOR:06}
D.~Bertsimas and A.~Thiele.
\newblock Robust and data-driven optimization: Modern decision making under
  uncertainty.
\newblock {\em Tutorials in Operations Research}, 4:95--122, 04 2006.

\bibitem{bertsimas:ieee:13}
Dimitris Bertsimas, Eugene Litvinov, Xu~Sun, Jinye Zhao, and Tongxin Zheng.
\newblock Adaptive robust optimization for the security constrained unit
  commitment problem.
\newblock {\em IEEE Transactions on Power Systems}, 28:52--63, 02 2013.

\bibitem{ElGhaoui:siam:97}
L.~El~Ghaoui and H.~Lebret.
\newblock Robust solutions to least-squares problems with uncertain data.
\newblock {\em SIAM Journal on Matrix Analysis and Applications},
  18(4):1035--1064, 1997.

\bibitem{Ghaoui:SIAM:98}
L.~El~Ghaoui, F.~Oustry, and H.~Lebret.
\newblock Robust solutions to uncertain semidefinite programs.
\newblock {\em SIAM Journal on Optimization}, 9(1):33--52, 1998.

\bibitem{gurobi:7}
{Gurobi Optimization, LLC}.
\newblock Gurobi optimizer reference manual, 2020.

\bibitem{Levin:OR:08}
Y.~Levin, J.~McGill, and M.~Nediak.
\newblock Risk in revenue management and dynamic pricing.
\newblock {\em Operations Research}, 56:326--343, 04 2008.

\bibitem{Mai:arxive:19}
T.~Mai and P.~Jaillet.
\newblock Robust multi-product pricing under general extreme value models,
  2019, 1912.09552.

\bibitem{mccormick:MP:75}
G.P. McCormick.
\newblock Computability of global solutions to factorable nonconvex programs:
  Part {I} — {C}onvex underestimating problems.
\newblock {\em Mathematical Programming}, 10:147--175, 1976.

\bibitem{Mulvey:OR:1995}
J.M. Mulvey, R.J. Vanderbei, and S.A. Zenios.
\newblock Robust optimization of large-scale systems.
\newblock {\em Operations Research}, 43(2):264--281, 1995.

\bibitem{MG:Nahmias:05}
S.~Nahmias.
\newblock {\em Production and Operations Analysis}.
\newblock McGraw-Hill, New York, 5th edition, 2005.

\bibitem{Ox:ozer:12}
{\" O}.~{\" O}zer and R.~Phillips.
\newblock {\em The Oxford Handbook of Pricing Management}.
\newblock Oxford University Press, 2012.

\bibitem{PRO:RoberPhillips:2005}
R.~Phillips.
\newblock {\em Pricing and Revenue Optimization}.
\newblock Stanford University Press, Stanford, CA, 2005.

\bibitem{MIT:Sheffi:05}
Y.~Sheffi.
\newblock {\em The Resilient Enterprise: Overcoming Vulnerability for
  Competitive Advantage}.
\newblock MIT Press, Cambridge, MA, 2005.

\bibitem{MG:Simchi:04}
D.~Simchi-Levi, P.~Kaminsky, and Simchi-Levi E.
\newblock {\em Managing the Supply Chain: The Definitive Guide for the Business
  Professional}.
\newblock McGraw-Hill, New York, 2004.

\bibitem{Soyster:OR:73}
A.L. Soyster.
\newblock Convex programming with set-inclusive constraints and applications to
  inexact linear programming.
\newblock {\em Operations Research}, 21(5):1154--1157, 1973.

\bibitem{Springer:talluri:05}
K.T. Talluri and G.J. van Ryzin.
\newblock {\em The Theory and Practice of Revenue Management}.
\newblock International Series in Operations Research \& Management Science.
  Springer US, Boston, MA, 2005.

\bibitem{thiele:OO:06}
A.~Thiele.
\newblock Single-product pricing via robust optimization.
\newblock Working Paper, 2006.

\bibitem{thiele:jrpm:09}
A.~Thiele.
\newblock Multi-product pricing via robust optimization.
\newblock {\em Journal of Revenue and Pricing Management}, 8:67--80, 2009.

\end{thebibliography}

\section*{Appendix}
\label{appx:appendx1}

Table~\ref{tbl:deltaoptval} presents each $\Delta$ value with its associated
optimal objective function value (revenue) of each instance. For ease of exposition, the optimal revenue value (denoted as $Z^{*}$) is normalized so the base scenario with $\Delta=0$ is set to zero. As we also expected, the larger the $\Delta$ value, the smaller $Z^{*}_{\Delta}$, as the problem becomes more conservative.

\begin{table} 
\begin{center}
	\footnotesize   
\resizebox{0.718\textwidth}{!}{\begin{tabular}{|l|c|rrrrrrrrrrr|r|}
	\hline\hline 
	& & \multicolumn{11}{|c|}{Optimal Values} & \\
	\cline{3-13}
	Sample & & 1 & 2 & 3 & 4 & 5 & 6 & 7 & 8 & 9 & 10 & 11 & Decline (\%) \\
\hline
\multirow{2}{*}{S1} & $\Delta $ & 0.0 & 0.65 & 1.3 & 1.9 & 2.6 & 3.2 & 3.9 & 4.5 & 5.2 & 5.8 & 6.5 & \multirow{2}{*}{4.0869} \\
& $Z^*$ & {\hspace{-2mm}\scriptsize 1.000000} & {\hspace{-2mm}\scriptsize 0.981036} & {\hspace{-2mm}\scriptsize 0.964493} & {\hspace{-2mm}\scriptsize 0.959536} & {\hspace{-2mm}\scriptsize 0.959145} & {\hspace{-2mm}\scriptsize 0.959132} & {\hspace{-2mm}\scriptsize 0.959131} & {\hspace{-2mm}\scriptsize 0.959131} & {\hspace{-2mm}\scriptsize 0.959131} & {\hspace{-2mm}\scriptsize 0.959131} & {\hspace{-2mm}\scriptsize 0.959131} & \\
\hline
\multirow{2}{*}{S2} & $\Delta $ & 0.0 & 0.55 & 1.1 & 1.7 & 2.2 & 2.8 & 3.3 & 3.9 & 4.4 & 5.0 & 5.5 & \multirow{2}{*}{0.7784} \\
& $Z^*$ & {\hspace{-2mm}\scriptsize 1.000000} & {\hspace{-2mm}\scriptsize 0.997144} & {\hspace{-2mm}\scriptsize 0.994409} & {\hspace{-2mm}\scriptsize 0.992639} & {\hspace{-2mm}\scriptsize 0.992311} & {\hspace{-2mm}\scriptsize 0.992245} & {\hspace{-2mm}\scriptsize 0.992220} & {\hspace{-2mm}\scriptsize 0.992216} & {\hspace{-2mm}\scriptsize 0.992216} & {\hspace{-2mm}\scriptsize 0.992216} & {\hspace{-2mm}\scriptsize 0.992216} & \\
\hline
\multirow{2}{*}{S3} & $\Delta $ & 0.0 & 0.46 & 0.92 & 1.4 & 1.8 & 2.3 & 2.8 & 3.2 & 3.7 & 4.1 & 4.6 & \multirow{2}{*}{6.1058} \\
& $Z^*$ & {\hspace{-2mm}\scriptsize 1.000000} & {\hspace{-2mm}\scriptsize 0.980973} & {\hspace{-2mm}\scriptsize 0.962698} & {\hspace{-2mm}\scriptsize 0.948051} & {\hspace{-2mm}\scriptsize 0.940493} & {\hspace{-2mm}\scriptsize 0.939001} & {\hspace{-2mm}\scriptsize 0.938953} & {\hspace{-2mm}\scriptsize 0.938943} & {\hspace{-2mm}\scriptsize 0.938942} & {\hspace{-2mm}\scriptsize 0.938942} & {\hspace{-2mm}\scriptsize 0.938942} & \\
\hline
\multirow{2}{*}{S4} & $\Delta $ & 0.0 & 0.39 & 0.78 & 1.2 & 1.6 & 1.9 & 2.3 & 2.7 & 3.1 & 3.5 & 3.9 & \multirow{2}{*}{5.2174} \\
& $Z^*$ & {\hspace{-2mm}\scriptsize 1.000000} & {\hspace{-2mm}\scriptsize 0.986918} & {\hspace{-2mm}\scriptsize 0.973921} & {\hspace{-2mm}\scriptsize 0.960932} & {\hspace{-2mm}\scriptsize 0.948030} & {\hspace{-2mm}\scriptsize 0.947895} & {\hspace{-2mm}\scriptsize 0.947826} & {\hspace{-2mm}\scriptsize 0.947826} & {\hspace{-2mm}\scriptsize 0.947826} & {\hspace{-2mm}\scriptsize 0.947826} & {\hspace{-2mm}\scriptsize 0.947826} & \\
\hline
\multirow{2}{*}{S5} & $\Delta $ & 0.0 & 0.35 & 0.71 & 1.1 & 1.4 & 1.8 & 2.1 & 2.5 & 2.8 & 3.2 & 3.5 & \multirow{2}{*}{0.0964} \\
& $Z^*$ & {\hspace{-2mm}\scriptsize 1.000000} & {\hspace{-2mm}\scriptsize 0.999791} & {\hspace{-2mm}\scriptsize 0.999583} & {\hspace{-2mm}\scriptsize 0.999374} & {\hspace{-2mm}\scriptsize 0.999165} & {\hspace{-2mm}\scriptsize 0.999040} & {\hspace{-2mm}\scriptsize 0.999036} & {\hspace{-2mm}\scriptsize 0.999036} & {\hspace{-2mm}\scriptsize 0.999036} & {\hspace{-2mm}\scriptsize 0.999036} & {\hspace{-2mm}\scriptsize 0.999036} & \\
\hline
\multirow{2}{*}{S6} & $\Delta $ & 0.0 & 0.4 & 0.8 & 1.2 & 1.6 & 2.0 & 2.4 & 2.8 & 3.2 & 3.6 & 4.0 & \multirow{2}{*}{0.2366} \\
& $Z^*$ & {\hspace{-2mm}\scriptsize 1.000000} & {\hspace{-2mm}\scriptsize 0.999404} & {\hspace{-2mm}\scriptsize 0.998809} & {\hspace{-2mm}\scriptsize 0.998223} & {\hspace{-2mm}\scriptsize 0.997753} & {\hspace{-2mm}\scriptsize 0.997637} & {\hspace{-2mm}\scriptsize 0.997634} & {\hspace{-2mm}\scriptsize 0.997634} & {\hspace{-2mm}\scriptsize 0.997634} & {\hspace{-2mm}\scriptsize 0.997634} & {\hspace{-2mm}\scriptsize 0.997634} & \\
\hline
\multirow{2}{*}{S7} & $\Delta $ & 0.0 & 0.4 & 0.8 & 1.2 & 1.6 & 2.0 & 2.4 & 2.8 & 3.2 & 3.6 & 4.0 & \multirow{2}{*}{0.0844} \\
& $Z^*$ & {\hspace{-2mm}\scriptsize 1.000000} & {\hspace{-2mm}\scriptsize 0.999798} & {\hspace{-2mm}\scriptsize 0.999596} & {\hspace{-2mm}\scriptsize 0.999395} & {\hspace{-2mm}\scriptsize 0.999204} & {\hspace{-2mm}\scriptsize 0.999156} & {\hspace{-2mm}\scriptsize 0.999156} & {\hspace{-2mm}\scriptsize 0.999156} & {\hspace{-2mm}\scriptsize 0.999156} & {\hspace{-2mm}\scriptsize 0.999156} & {\hspace{-2mm}\scriptsize 0.999156} & \\
\hline
\multirow{2}{*}{S8} & $\Delta $ & 0.0 & 0.39 & 0.77 & 1.2 & 1.5 & 1.9 & 2.3 & 2.7 & 3.1 & 3.5 & 3.9 & \multirow{2}{*}{0.1391} \\
& $Z^*$ & {\hspace{-2mm}\scriptsize 1.000000} & {\hspace{-2mm}\scriptsize 0.999667} & {\hspace{-2mm}\scriptsize 0.999334} & {\hspace{-2mm}\scriptsize 0.999002} & {\hspace{-2mm}\scriptsize 0.998704} & {\hspace{-2mm}\scriptsize 0.998611} & {\hspace{-2mm}\scriptsize 0.998609} & {\hspace{-2mm}\scriptsize 0.998609} & {\hspace{-2mm}\scriptsize 0.998609} & {\hspace{-2mm}\scriptsize 0.998609} & {\hspace{-2mm}\scriptsize 0.998609} & \\
\hline
\multirow{2}{*}{S9} & $\Delta $ & 0.0 & 0.4 & 0.8 & 1.2 & 1.6 & 2.0 & 2.4 & 2.8 & 3.2 & 3.6 & 4.0 & \multirow{2}{*}{3.5418} \\
& $Z^*$ & {\hspace{-2mm}\scriptsize 1.000000} & {\hspace{-2mm}\scriptsize 0.988791} & {\hspace{-2mm}\scriptsize 0.977581} & {\hspace{-2mm}\scriptsize 0.966854} & {\hspace{-2mm}\scriptsize 0.964664} & {\hspace{-2mm}\scriptsize 0.964583} & {\hspace{-2mm}\scriptsize 0.964582} & {\hspace{-2mm}\scriptsize 0.964582} & {\hspace{-2mm}\scriptsize 0.964582} & {\hspace{-2mm}\scriptsize 0.964582} & {\hspace{-2mm}\scriptsize 0.964582} & \\
\hline
\multirow{2}{*}{S10} & $\Delta $ & 0.0 & 0.42 & 0.83 & 1.2 & 1.7 & 2.1 & 2.5 & 2.9 & 3.3 & 3.7 & 4.2 & \multirow{2}{*}{0.1005} \\
& $Z^*$ & {\hspace{-2mm}\scriptsize 1.000000} & {\hspace{-2mm}\scriptsize 0.999753} & {\hspace{-2mm}\scriptsize 0.999507} & {\hspace{-2mm}\scriptsize 0.999260} & {\hspace{-2mm}\scriptsize 0.999055} & {\hspace{-2mm}\scriptsize 0.998999} & {\hspace{-2mm}\scriptsize 0.998995} & {\hspace{-2mm}\scriptsize 0.998995} & {\hspace{-2mm}\scriptsize 0.998995} & {\hspace{-2mm}\scriptsize 0.998995} & {\hspace{-2mm}\scriptsize 0.998995} & \\
\hline
\multirow{2}{*}{S11} & $\Delta $ & 0.0 & 0.39 & 0.77 & 1.2 & 1.5 & 1.9 & 2.3 & 2.7 & 3.1 & 3.5 & 3.9 & \multirow{2}{*}{0.3744} \\
& $Z^*$ & {\hspace{-2mm}\scriptsize 1.000000} & {\hspace{-2mm}\scriptsize 0.999164} & {\hspace{-2mm}\scriptsize 0.998361} & {\hspace{-2mm}\scriptsize 0.997558} & {\hspace{-2mm}\scriptsize 0.996762} & {\hspace{-2mm}\scriptsize 0.996259} & {\hspace{-2mm}\scriptsize 0.996256} & {\hspace{-2mm}\scriptsize 0.996256} & {\hspace{-2mm}\scriptsize 0.996256} & {\hspace{-2mm}\scriptsize 0.996256} & {\hspace{-2mm}\scriptsize 0.996256} & \\
\hline
\multirow{2}{*}{S12} & $\Delta $ & 0.0 & 0.39 & 0.77 & 1.2 & 1.5 & 1.9 & 2.3 & 2.7 & 3.1 & 3.5 & 3.9 & \multirow{2}{*}{0.1979} \\
& $Z^*$ & {\hspace{-2mm}\scriptsize 1.000000} & {\hspace{-2mm}\scriptsize 0.999501} & {\hspace{-2mm}\scriptsize 0.999001} & {\hspace{-2mm}\scriptsize 0.998502} & {\hspace{-2mm}\scriptsize 0.998129} & {\hspace{-2mm}\scriptsize 0.998021} & {\hspace{-2mm}\scriptsize 0.998021} & {\hspace{-2mm}\scriptsize 0.998021} & {\hspace{-2mm}\scriptsize 0.998021} & {\hspace{-2mm}\scriptsize 0.998021} & {\hspace{-2mm}\scriptsize 0.998021} & \\
\hline
\multirow{2}{*}{S13} & $\Delta $ & 0.0 & 0.41 & 0.83 & 1.2 & 1.7 & 2.1 & 2.5 & 2.9 & 3.3 & 3.7 & 4.1 & \multirow{2}{*}{0.1200} \\
& $Z^*$ & {\hspace{-2mm}\scriptsize 1.000000} & {\hspace{-2mm}\scriptsize 0.999697} & {\hspace{-2mm}\scriptsize 0.999394} & {\hspace{-2mm}\scriptsize 0.999091} & {\hspace{-2mm}\scriptsize 0.998829} & {\hspace{-2mm}\scriptsize 0.998801} & {\hspace{-2mm}\scriptsize 0.998800} & {\hspace{-2mm}\scriptsize 0.998800} & {\hspace{-2mm}\scriptsize 0.998800} & {\hspace{-2mm}\scriptsize 0.998800} & {\hspace{-2mm}\scriptsize 0.998800} & \\
\hline
\multirow{2}{*}{S14} & $\Delta $ & 0.0 & 0.38 & 0.77 & 1.1 & 1.5 & 1.9 & 2.3 & 2.7 & 3.1 & 3.4 & 3.8 & \multirow{2}{*}{0.1398} \\
& $Z^*$ & {\hspace{-2mm}\scriptsize 1.000000} & {\hspace{-2mm}\scriptsize 0.999686} & {\hspace{-2mm}\scriptsize 0.999372} & {\hspace{-2mm}\scriptsize 0.999058} & {\hspace{-2mm}\scriptsize 0.998750} & {\hspace{-2mm}\scriptsize 0.998604} & {\hspace{-2mm}\scriptsize 0.998602} & {\hspace{-2mm}\scriptsize 0.998602} & {\hspace{-2mm}\scriptsize 0.998602} & {\hspace{-2mm}\scriptsize 0.998602} & {\hspace{-2mm}\scriptsize 0.998602} & \\
\hline
\multirow{2}{*}{S15} & $\Delta $ & 0.0 & 0.48 & 0.95 & 1.4 & 1.9 & 2.4 & 2.9 & 3.3 & 3.8 & 4.3 & 4.8 & \multirow{2}{*}{0.3392} \\
& $Z^*$ & {\hspace{-2mm}\scriptsize 1.000000} & {\hspace{-2mm}\scriptsize 0.998885} & {\hspace{-2mm}\scriptsize 0.997798} & {\hspace{-2mm}\scriptsize 0.996985} & {\hspace{-2mm}\scriptsize 0.996681} & {\hspace{-2mm}\scriptsize 0.996632} & {\hspace{-2mm}\scriptsize 0.996617} & {\hspace{-2mm}\scriptsize 0.996609} & {\hspace{-2mm}\scriptsize 0.996609} & {\hspace{-2mm}\scriptsize 0.996608} & {\hspace{-2mm}\scriptsize 0.996608} & \\
\hline
\multirow{2}{*}{S16} & $\Delta $ & 0.0 & 0.41 & 0.83 & 1.2 & 1.7 & 2.1 & 2.5 & 2.9 & 3.3 & 3.7 & 4.1 & \multirow{2}{*}{0.0951} \\
& $Z^*$ & {\hspace{-2mm}\scriptsize 1.000000} & {\hspace{-2mm}\scriptsize 0.999764} & {\hspace{-2mm}\scriptsize 0.999529} & {\hspace{-2mm}\scriptsize 0.999293} & {\hspace{-2mm}\scriptsize 0.999094} & {\hspace{-2mm}\scriptsize 0.999049} & {\hspace{-2mm}\scriptsize 0.999049} & {\hspace{-2mm}\scriptsize 0.999049} & {\hspace{-2mm}\scriptsize 0.999049} & {\hspace{-2mm}\scriptsize 0.999049} & {\hspace{-2mm}\scriptsize 0.999049} & \\
\hline
\multirow{2}{*}{S17} & $\Delta $ & 0.0 & 0.41 & 0.83 & 1.2 & 1.7 & 2.1 & 2.5 & 2.9 & 3.3 & 3.7 & 4.1 & \multirow{2}{*}{0.1019} \\
& $Z^*$ & {\hspace{-2mm}\scriptsize 1.000000} & {\hspace{-2mm}\scriptsize 0.999742} & {\hspace{-2mm}\scriptsize 0.999484} & {\hspace{-2mm}\scriptsize 0.999226} & {\hspace{-2mm}\scriptsize 0.999010} & {\hspace{-2mm}\scriptsize 0.998982} & {\hspace{-2mm}\scriptsize 0.998981} & {\hspace{-2mm}\scriptsize 0.998981} & {\hspace{-2mm}\scriptsize 0.998981} & {\hspace{-2mm}\scriptsize 0.998981} & {\hspace{-2mm}\scriptsize 0.998981} & \\
\hline
\multirow{2}{*}{S18} & $\Delta $ & 0.0 & 0.41 & 0.82 & 1.2 & 1.6 & 2.1 & 2.5 & 2.9 & 3.3 & 3.7 & 4.1 & \multirow{2}{*}{0.0716} \\
& $Z^*$ & {\hspace{-2mm}\scriptsize 1.000000} & {\hspace{-2mm}\scriptsize 0.999821} & {\hspace{-2mm}\scriptsize 0.999642} & {\hspace{-2mm}\scriptsize 0.999464} & {\hspace{-2mm}\scriptsize 0.999313} & {\hspace{-2mm}\scriptsize 0.999285} & {\hspace{-2mm}\scriptsize 0.999284} & {\hspace{-2mm}\scriptsize 0.999284} & {\hspace{-2mm}\scriptsize 0.999284} & {\hspace{-2mm}\scriptsize 0.999284} & {\hspace{-2mm}\scriptsize 0.999284} & \\
\hline
\multirow{2}{*}{S19} & $\Delta $ & 0.0 & 0.44 & 0.88 & 1.3 & 1.8 & 2.2 & 2.6 & 3.1 & 3.5 & 4.0 & 4.4 & \multirow{2}{*}{0.1092} \\
& $Z^*$ & {\hspace{-2mm}\scriptsize 1.000000} & {\hspace{-2mm}\scriptsize 0.999714} & {\hspace{-2mm}\scriptsize 0.999429} & {\hspace{-2mm}\scriptsize 0.999143} & {\hspace{-2mm}\scriptsize 0.998933} & {\hspace{-2mm}\scriptsize 0.998908} & {\hspace{-2mm}\scriptsize 0.998908} & {\hspace{-2mm}\scriptsize 0.998908} & {\hspace{-2mm}\scriptsize 0.998908} & {\hspace{-2mm}\scriptsize 0.998908} & {\hspace{-2mm}\scriptsize 0.998908} & \\
\hline
\multirow{2}{*}{S20} & $\Delta $ & 0.0 & 0.38 & 0.77 & 1.2 & 1.5 & 1.9 & 2.3 & 2.7 & 3.1 & 3.5 & 3.8 & \multirow{2}{*}{0.1906} \\
& $Z^*$ & {\hspace{-2mm}\scriptsize 1.000000} & {\hspace{-2mm}\scriptsize 0.999560} & {\hspace{-2mm}\scriptsize 0.999120} & {\hspace{-2mm}\scriptsize 0.998679} & {\hspace{-2mm}\scriptsize 0.998257} & {\hspace{-2mm}\scriptsize 0.998097} & {\hspace{-2mm}\scriptsize 0.998094} & {\hspace{-2mm}\scriptsize 0.998094} & {\hspace{-2mm}\scriptsize 0.998094} & {\hspace{-2mm}\scriptsize 0.998094} & {\hspace{-2mm}\scriptsize 0.998094} & \\
\hline
\multirow{2}{*}{S21} & $\Delta $ & 0.0 & 0.39 & 0.77 & 1.2 & 1.5 & 1.9 & 2.3 & 2.7 & 3.1 & 3.5 & 3.9 & \multirow{2}{*}{0.0898} \\
& $Z^*$ & {\hspace{-2mm}\scriptsize 1.000000} & {\hspace{-2mm}\scriptsize 0.999793} & {\hspace{-2mm}\scriptsize 0.999585} & {\hspace{-2mm}\scriptsize 0.999378} & {\hspace{-2mm}\scriptsize 0.999176} & {\hspace{-2mm}\scriptsize 0.999104} & {\hspace{-2mm}\scriptsize 0.999102} & {\hspace{-2mm}\scriptsize 0.999102} & {\hspace{-2mm}\scriptsize 0.999102} & {\hspace{-2mm}\scriptsize 0.999102} & {\hspace{-2mm}\scriptsize 0.999102} & \\
\hline
\multirow{2}{*}{S22} & $\Delta $ & 0.0 & 0.41 & 0.81 & 1.2 & 1.6 & 2.0 & 2.4 & 2.8 & 3.3 & 3.7 & 4.1 & \multirow{2}{*}{0.0757} \\
& $Z^*$ & {\hspace{-2mm}\scriptsize 1.000000} & {\hspace{-2mm}\scriptsize 0.999818} & {\hspace{-2mm}\scriptsize 0.999635} & {\hspace{-2mm}\scriptsize 0.999453} & {\hspace{-2mm}\scriptsize 0.999292} & {\hspace{-2mm}\scriptsize 0.999246} & {\hspace{-2mm}\scriptsize 0.999243} & {\hspace{-2mm}\scriptsize 0.999243} & {\hspace{-2mm}\scriptsize 0.999243} & {\hspace{-2mm}\scriptsize 0.999243} & {\hspace{-2mm}\scriptsize 0.999243} & \\
\hline
\multirow{2}{*}{S23} & $\Delta $ & 0.0 & 0.39 & 0.79 & 1.2 & 1.6 & 2.0 & 2.4 & 2.8 & 3.2 & 3.6 & 3.9 & \multirow{2}{*}{0.2508} \\
& $Z^*$ & {\hspace{-2mm}\scriptsize 1.000000} & {\hspace{-2mm}\scriptsize 0.999300} & {\hspace{-2mm}\scriptsize 0.998601} & {\hspace{-2mm}\scriptsize 0.997901} & {\hspace{-2mm}\scriptsize 0.997588} & {\hspace{-2mm}\scriptsize 0.997496} & {\hspace{-2mm}\scriptsize 0.997492} & {\hspace{-2mm}\scriptsize 0.997492} & {\hspace{-2mm}\scriptsize 0.997492} & {\hspace{-2mm}\scriptsize 0.997492} & {\hspace{-2mm}\scriptsize 0.997492} & \\
\hline
\multirow{2}{*}{S24} & $\Delta $ & 0.0 & 0.41 & 0.82 & 1.2 & 1.6 & 2.0 & 2.5 & 2.9 & 3.3 & 3.7 & 4.1 & \multirow{2}{*}{0.0681} \\
& $Z^*$ & {\hspace{-2mm}\scriptsize 1.000000} & {\hspace{-2mm}\scriptsize 0.999831} & {\hspace{-2mm}\scriptsize 0.999663} & {\hspace{-2mm}\scriptsize 0.999494} & {\hspace{-2mm}\scriptsize 0.999352} & {\hspace{-2mm}\scriptsize 0.999320} & {\hspace{-2mm}\scriptsize 0.999319} & {\hspace{-2mm}\scriptsize 0.999319} & {\hspace{-2mm}\scriptsize 0.999319} & {\hspace{-2mm}\scriptsize 0.999319} & {\hspace{-2mm}\scriptsize 0.999319} & \\
\hline
\multirow{2}{*}{S25} & $\Delta $ & 0.0 & 0.62 & 1.2 & 1.9 & 2.5 & 3.1 & 3.7 & 4.3 & 4.9 & 5.6 & 6.2 & \multirow{2}{*}{0.8945} \\
& $Z^*$ & {\hspace{-2mm}\scriptsize 1.000000} & {\hspace{-2mm}\scriptsize 0.996451} & {\hspace{-2mm}\scriptsize 0.992945} & {\hspace{-2mm}\scriptsize 0.991132} & {\hspace{-2mm}\scriptsize 0.991061} & {\hspace{-2mm}\scriptsize 0.991056} & {\hspace{-2mm}\scriptsize 0.991055} & {\hspace{-2mm}\scriptsize 0.991055} & {\hspace{-2mm}\scriptsize 0.991055} & {\hspace{-2mm}\scriptsize 0.991055} & {\hspace{-2mm}\scriptsize 0.991055} & \\
\hline
\multirow{2}{*}{S26} & $\Delta $ & 0.0 & 0.65 & 1.3 & 1.9 & 2.6 & 3.2 & 3.9 & 4.5 & 5.2 & 5.8 & 6.5 & \multirow{2}{*}{1.8451} \\
& $Z^*$ & {\hspace{-2mm}\scriptsize 1.000000} & {\hspace{-2mm}\scriptsize 0.992429} & {\hspace{-2mm}\scriptsize 0.985244} & {\hspace{-2mm}\scriptsize 0.982216} & {\hspace{-2mm}\scriptsize 0.981786} & {\hspace{-2mm}\scriptsize 0.981554} & {\hspace{-2mm}\scriptsize 0.981550} & {\hspace{-2mm}\scriptsize 0.981550} & {\hspace{-2mm}\scriptsize 0.981549} & {\hspace{-2mm}\scriptsize 0.981549} & {\hspace{-2mm}\scriptsize 0.981549} & \\
\hline
\multirow{2}{*}{S27} & $\Delta $ & 0.0 & 0.57 & 1.1 & 1.7 & 2.3 & 2.9 & 3.4 & 4.0 & 4.6 & 5.2 & 5.7 & \multirow{2}{*}{0.0501} \\
& $Z^*$ & {\hspace{-2mm}\scriptsize 1.000000} & {\hspace{-2mm}\scriptsize 0.999812} & {\hspace{-2mm}\scriptsize 0.999630} & {\hspace{-2mm}\scriptsize 0.999518} & {\hspace{-2mm}\scriptsize 0.999504} & {\hspace{-2mm}\scriptsize 0.999501} & {\hspace{-2mm}\scriptsize 0.999500} & {\hspace{-2mm}\scriptsize 0.999500} & {\hspace{-2mm}\scriptsize 0.999500} & {\hspace{-2mm}\scriptsize 0.999499} & {\hspace{-2mm}\scriptsize 0.999499} & \\
\hline
\multirow{2}{*}{S28} & $\Delta $ & 0.0 & 0.62 & 1.2 & 1.9 & 2.5 & 3.1 & 3.7 & 4.4 & 5.0 & 5.6 & 6.2 & \multirow{2}{*}{0.0691} \\
& $Z^*$ & {\hspace{-2mm}\scriptsize 1.000000} & {\hspace{-2mm}\scriptsize 0.999720} & {\hspace{-2mm}\scriptsize 0.999459} & {\hspace{-2mm}\scriptsize 0.999331} & {\hspace{-2mm}\scriptsize 0.999316} & {\hspace{-2mm}\scriptsize 0.999311} & {\hspace{-2mm}\scriptsize 0.999309} & {\hspace{-2mm}\scriptsize 0.999309} & {\hspace{-2mm}\scriptsize 0.999309} & {\hspace{-2mm}\scriptsize 0.999309} & {\hspace{-2mm}\scriptsize 0.999309} & \\
\hline
\multirow{2}{*}{S29} & $\Delta $ & 0.0 & 0.62 & 1.2 & 1.8 & 2.5 & 3.1 & 3.7 & 4.3 & 4.9 & 5.5 & 6.2 & \multirow{2}{*}{2.4146} \\
& $Z^*$ & {\hspace{-2mm}\scriptsize 1.000000} & {\hspace{-2mm}\scriptsize 0.990004} & {\hspace{-2mm}\scriptsize 0.981477} & {\hspace{-2mm}\scriptsize 0.977098} & {\hspace{-2mm}\scriptsize 0.976042} & {\hspace{-2mm}\scriptsize 0.975912} & {\hspace{-2mm}\scriptsize 0.975857} & {\hspace{-2mm}\scriptsize 0.975854} & {\hspace{-2mm}\scriptsize 0.975854} & {\hspace{-2mm}\scriptsize 0.975854} & {\hspace{-2mm}\scriptsize 0.975854} & \\
\hline
\multirow{2}{*}{S30} & $\Delta $ & 0.0 & 0.6 & 1.2 & 1.8 & 2.4 & 3.0 & 3.6 & 4.2 & 4.8 & 5.4 & 6.0 & \multirow{2}{*}{0.5939} \\
& $Z^*$ & {\hspace{-2mm}\scriptsize 1.000000} & {\hspace{-2mm}\scriptsize 0.997535} & {\hspace{-2mm}\scriptsize 0.995168} & {\hspace{-2mm}\scriptsize 0.994144} & {\hspace{-2mm}\scriptsize 0.994083} & {\hspace{-2mm}\scriptsize 0.994070} & {\hspace{-2mm}\scriptsize 0.994063} & {\hspace{-2mm}\scriptsize 0.994061} & {\hspace{-2mm}\scriptsize 0.994061} & {\hspace{-2mm}\scriptsize 0.994061} & {\hspace{-2mm}\scriptsize 0.994061} & \\
\hline
\multirow{2}{*}{S31} & $\Delta $ & 0.0 & 0.7 & 1.4 & 2.1 & 2.8 & 3.5 & 4.2 & 4.9 & 5.6 & 6.3 & 7.0 & \multirow{2}{*}{6.9645} \\
& $Z^*$ & {\hspace{-2mm}\scriptsize 1.000000} & {\hspace{-2mm}\scriptsize 0.968040} & {\hspace{-2mm}\scriptsize 0.940502} & {\hspace{-2mm}\scriptsize 0.930946} & {\hspace{-2mm}\scriptsize 0.930493} & {\hspace{-2mm}\scriptsize 0.930380} & {\hspace{-2mm}\scriptsize 0.930359} & {\hspace{-2mm}\scriptsize 0.930356} & {\hspace{-2mm}\scriptsize 0.930355} & {\hspace{-2mm}\scriptsize 0.930355} & {\hspace{-2mm}\scriptsize 0.930355} & \\
\hline
\multirow{2}{*}{S32} & $\Delta $ & 0.0 & 0.59 & 1.2 & 1.8 & 2.4 & 3.0 & 3.6 & 4.2 & 4.8 & 5.4 & 5.9 & \multirow{2}{*}{2.6059} \\
& $Z^*$ & {\hspace{-2mm}\scriptsize 1.000000} & {\hspace{-2mm}\scriptsize 0.989968} & {\hspace{-2mm}\scriptsize 0.981014} & {\hspace{-2mm}\scriptsize 0.975376} & {\hspace{-2mm}\scriptsize 0.974279} & {\hspace{-2mm}\scriptsize 0.974051} & {\hspace{-2mm}\scriptsize 0.973967} & {\hspace{-2mm}\scriptsize 0.973953} & {\hspace{-2mm}\scriptsize 0.973945} & {\hspace{-2mm}\scriptsize 0.973941} & {\hspace{-2mm}\scriptsize 0.973941} & \\
\hline
\multirow{2}{*}{S33} & $\Delta $ & 0.0 & 0.58 & 1.2 & 1.7 & 2.3 & 2.9 & 3.5 & 4.1 & 4.7 & 5.2 & 5.8 & \multirow{2}{*}{1.2579} \\
& $Z^*$ & {\hspace{-2mm}\scriptsize 1.000000} & {\hspace{-2mm}\scriptsize 0.995042} & {\hspace{-2mm}\scriptsize 0.990318} & {\hspace{-2mm}\scriptsize 0.987787} & {\hspace{-2mm}\scriptsize 0.987447} & {\hspace{-2mm}\scriptsize 0.987430} & {\hspace{-2mm}\scriptsize 0.987423} & {\hspace{-2mm}\scriptsize 0.987421} & {\hspace{-2mm}\scriptsize 0.987421} & {\hspace{-2mm}\scriptsize 0.987421} & {\hspace{-2mm}\scriptsize 0.987421} & \\
\hline
\multirow{2}{*}{S34} & $\Delta $ & 0.0 & 0.57 & 1.1 & 1.7 & 2.3 & 2.8 & 3.4 & 4.0 & 4.5 & 5.1 & 5.7 & \multirow{2}{*}{0.2168} \\
& $Z^*$ & {\hspace{-2mm}\scriptsize 1.000000} & {\hspace{-2mm}\scriptsize 0.999182} & {\hspace{-2mm}\scriptsize 0.998445} & {\hspace{-2mm}\scriptsize 0.997991} & {\hspace{-2mm}\scriptsize 0.997875} & {\hspace{-2mm}\scriptsize 0.997843} & {\hspace{-2mm}\scriptsize 0.997833} & {\hspace{-2mm}\scriptsize 0.997832} & {\hspace{-2mm}\scriptsize 0.997832} & {\hspace{-2mm}\scriptsize 0.997832} & {\hspace{-2mm}\scriptsize 0.997832} & \\
\hline
\multirow{2}{*}{S35} & $\Delta $ & 0.0 & 0.65 & 1.3 & 1.9 & 2.6 & 3.2 & 3.9 & 4.5 & 5.2 & 5.8 & 6.5 & \multirow{2}{*}{0.1478} \\
& $Z^*$ & {\hspace{-2mm}\scriptsize 1.000000} & {\hspace{-2mm}\scriptsize 0.999383} & {\hspace{-2mm}\scriptsize 0.998828} & {\hspace{-2mm}\scriptsize 0.998545} & {\hspace{-2mm}\scriptsize 0.998530} & {\hspace{-2mm}\scriptsize 0.998523} & {\hspace{-2mm}\scriptsize 0.998522} & {\hspace{-2mm}\scriptsize 0.998522} & {\hspace{-2mm}\scriptsize 0.998522} & {\hspace{-2mm}\scriptsize 0.998522} & {\hspace{-2mm}\scriptsize 0.998522} & \\
\hline
\multirow{2}{*}{S36} & $\Delta $ & 0.0 & 0.56 & 1.1 & 1.7 & 2.2 & 2.8 & 3.4 & 3.9 & 4.5 & 5.1 & 5.6 & \multirow{2}{*}{0.8719} \\
& $Z^*$ & {\hspace{-2mm}\scriptsize 1.000000} & {\hspace{-2mm}\scriptsize 0.996699} & {\hspace{-2mm}\scriptsize 0.993594} & {\hspace{-2mm}\scriptsize 0.991802} & {\hspace{-2mm}\scriptsize 0.991374} & {\hspace{-2mm}\scriptsize 0.991324} & {\hspace{-2mm}\scriptsize 0.991304} & {\hspace{-2mm}\scriptsize 0.991291} & {\hspace{-2mm}\scriptsize 0.991284} & {\hspace{-2mm}\scriptsize 0.991281} & {\hspace{-2mm}\scriptsize 0.991281} & \\
\hline
\multirow{2}{*}{S37} & $\Delta $ & 0.0 & 0.65 & 1.3 & 1.9 & 2.6 & 3.2 & 3.9 & 4.5 & 5.2 & 5.8 & 6.5 & \multirow{2}{*}{3.0746} \\
& $Z^*$ & {\hspace{-2mm}\scriptsize 1.000000} & {\hspace{-2mm}\scriptsize 0.985972} & {\hspace{-2mm}\scriptsize 0.972225} & {\hspace{-2mm}\scriptsize 0.969320} & {\hspace{-2mm}\scriptsize 0.969271} & {\hspace{-2mm}\scriptsize 0.969257} & {\hspace{-2mm}\scriptsize 0.969254} & {\hspace{-2mm}\scriptsize 0.969254} & {\hspace{-2mm}\scriptsize 0.969254} & {\hspace{-2mm}\scriptsize 0.969254} & {\hspace{-2mm}\scriptsize 0.969254} & \\
\hline
\multirow{2}{*}{S38} & $\Delta $ & 0.0 & 0.65 & 1.3 & 2.0 & 2.6 & 3.3 & 3.9 & 4.6 & 5.2 & 5.9 & 6.5 & \multirow{2}{*}{1.1271} \\
& $Z^*$ & {\hspace{-2mm}\scriptsize 1.000000} & {\hspace{-2mm}\scriptsize 0.995359} & {\hspace{-2mm}\scriptsize 0.990956} & {\hspace{-2mm}\scriptsize 0.988909} & {\hspace{-2mm}\scriptsize 0.988752} & {\hspace{-2mm}\scriptsize 0.988731} & {\hspace{-2mm}\scriptsize 0.988729} & {\hspace{-2mm}\scriptsize 0.988729} & {\hspace{-2mm}\scriptsize 0.988729} & {\hspace{-2mm}\scriptsize 0.988729} & {\hspace{-2mm}\scriptsize 0.988729} & \\
\hline
\multirow{2}{*}{S39} & $\Delta $ & 0.0 & 0.61 & 1.2 & 1.8 & 2.4 & 3.0 & 3.6 & 4.2 & 4.9 & 5.5 & 6.1 & \multirow{2}{*}{0.4145} \\
& $Z^*$ & {\hspace{-2mm}\scriptsize 1.000000} & {\hspace{-2mm}\scriptsize 0.998126} & {\hspace{-2mm}\scriptsize 0.996349} & {\hspace{-2mm}\scriptsize 0.995896} & {\hspace{-2mm}\scriptsize 0.995861} & {\hspace{-2mm}\scriptsize 0.995856} & {\hspace{-2mm}\scriptsize 0.995856} & {\hspace{-2mm}\scriptsize 0.995855} & {\hspace{-2mm}\scriptsize 0.995855} & {\hspace{-2mm}\scriptsize 0.995855} & {\hspace{-2mm}\scriptsize 0.995855} & \\
\hline
\multirow{2}{*}{S40} & $\Delta $ & 0.0 & 0.65 & 1.3 & 2.0 & 2.6 & 3.3 & 3.9 & 4.6 & 5.2 & 5.9 & 6.5 & \multirow{2}{*}{0.3026} \\
& $Z^*$ & {\hspace{-2mm}\scriptsize 1.000000} & {\hspace{-2mm}\scriptsize 0.998797} & {\hspace{-2mm}\scriptsize 0.997742} & {\hspace{-2mm}\scriptsize 0.997134} & {\hspace{-2mm}\scriptsize 0.996982} & {\hspace{-2mm}\scriptsize 0.996976} & {\hspace{-2mm}\scriptsize 0.996975} & {\hspace{-2mm}\scriptsize 0.996974} & {\hspace{-2mm}\scriptsize 0.996974} & {\hspace{-2mm}\scriptsize 0.996974} & {\hspace{-2mm}\scriptsize 0.996974} & \\
	\hline \hline 

\end{tabular}}
\caption{\small $\Delta$ and optimal values of the price optimization problem with robust approach. \label{tbl:deltaoptval}}
\end{center}
\end{table}

\end{document}